\documentstyle[12pt]{article}

\input amssym.def
\input amssym.tex

    \newtheorem{rema}{Remark}[section]
    \newtheorem{propo}[rema]{Proposition}
    \newtheorem{theo}[rema]{Theorem}
   \newtheorem{defi}[rema]{Definition}

     \newtheorem{exam}[rema]{Example}
        \newcommand{\ba}{\begin{array}}
        \newcommand{\ea}{\end{array}}
        \newcommand{\be}{\begin{equation}}
        \newcommand{\ee}{\end{equation}}
        \newcommand{\bea}{\begin{eqnarray}}
        \newcommand{\eea}{\end{eqnarray}}
        \newcommand{\nno}{\nonumber}
        
\newcommand{\mbar}{\mbox{\huge $\vert$}}
        \newcommand{\p}{\partial}
        \newcommand{\dps}{\displaystyle}

 \newcommand{\res}{\mbox{\rm Res}}
 
 \newcommand{\pf}{{\it Proof}\hspace{2ex}}
 \newcommand{\epf}{\hspace{2em}$\Box$}
 \newcommand{\epfv}{\hspace{1em}$\Box$\vspace{1em}}
\renewcommand{\hom}{\mbox{\rm Hom}}

    \title{{\bf On the ${\cal D}$-module and formal-variable 
approaches to vertex algebras}}
    \author{Yi-Zhi Huang and James Lepowsky}
    \date{}
    \begin{document}
    \bibliographystyle{alpha}
    \maketitle

\renewcommand{\theequation}{\thesection.\arabic{equation}}
\renewcommand{\therema}{\thesection.\arabic{rema}}
\setcounter{equation}{0}
\setcounter{rema}{0}

\section{Introduction}

In a program to formulate and develop two-dimensional conformal field
theory in the framework of algebraic geometry, Beilinson and Drinfeld
\cite{BD} have recently given a notion of ``chiral algebra'' in terms
of ${\cal D}$-modules on algebraic curves.  This definition consists
of a ``skew-symmetry'' relation and a ``Jacobi identity'' relation in
a categorical setting, and it leads to the operator product expansion
for holomorphic quantum fields in the spirit of two-dimensional
conformal field theory, as expressed in \cite{BPZ}.  Because this
operator product expansion, properly formulated, is known to be
essentially a variant of the main axiom, the ``Jacobi identity''
\cite{FLM}, for vertex (operator) algebras
(\cite{B}, \cite{FLM}; see \cite{FLM} for the proof), the chiral
algebras of \cite{BD} amount essentially to vertex algebras.

In this paper, we show directly that the chiral algebras of \cite{BD}
are essentially the same as vertex algebras without vacuum vector (and
without grading), by establishing an equivalence between the
skew-symmetry and Jacobi identity relations of \cite{BD} and the
(similarly-named, but different) skew-symmetry and Jacobi identity
relations in the formal-variable approach to vertex operator algebra
theory (see \cite{FLM}, \cite{FHL}).  In particular, among the equivalent
formulations of the notion of vertex (operator) algebra, the ${\cal
D}$-module notion of chiral algebra corresponds the most closely to
the formal-variable notion, rather than to, say, the
operator-product-expansion notion (based on the ``commutativity'' and
``associativity'' relations, as explained in \cite{FLM}, \cite{FHL}) or 
to the
geometric or operadic notion (\cite{H1}, \cite{H2}, \cite{HL}).

More precisely, we prove that for
any nonempty open subset $X$ of ${\Bbb
C}$, the category of vertex algebras without
vacuum over $X$ (see Definitions \ref{va-no-va} and
\ref{va-no-va-x}) and the category of chiral
algebras over $X$  (see Definition \ref{chiral})
are equivalent (see Section 5).
Beilinson-Drinfeld's notion of chiral algebra is in general formulated
over higher-genus curves.  But since chiral algebras in this sense are
essentially local objects, the equivalence proved in this paper shows
that the notion in \cite{BD} is indeed essentially equivalent to the
notion of vertex algebra without vacuum.

We hope that the present expository exercise helps to illuminate the
relations between the theories and philosophies of ${\cal D}$-modules
and of vertex operator algebras.  For example, in the Jacobi identity
for vertex algebras, the three formal variables are on equal footing
because of an intrinsic $S_3$-symmetry (see \cite{FHL}, Section 2.7),
while in the ${\cal D}$-module approach, there are only two (complex)
variables, as in the operator-product-expansion approach.  (See Remark
\ref{brokens3}.)  To see the $S_{3}$-symmetry explicitly in the
algebro-geometric framework, we would have to introduce an analogue of
the notion of ${\cal D}$-module allowing global translations of a
variable rather than just ``infinitesimal translations.''  The Jacobi
identity would then be interpreted as an identity in terms of such
``modified ${\cal D}$-modules,'' so that the three variables involved
would play symmetric roles.  This will be discussed in future
publications.  The Jacobi identity for vertex operator algebras and
its $S_{3}$-symmetry in fact play a central role in the theory of
vertex operator algebras, in particular, in the construction of
``vertex tensor categories'' (see \cite{HL4}, \cite{HL2}, \cite{HL3}).

Even without the introduction of such global translations of
variables, all of the many calculations in this paper involving
binomial expansions can be greatly simplified if we systematically
introduce formal (not complex) variables playing the role of ``formal
global translations.'' For instance, the expression
$(z_{1}-z_{2})^{n}(A_{1}\otimes A_{2})$, $n\in {\Bbb Z}$, occurring
starting in Section 4 can be viewed as the coefficient of $x^{-n-1}$
in $x^{-1}\delta(\frac{z_{1}-z_{2}}{x})(A_{1}\otimes A_{2})$, where
$x$ is a formal variable and
$\delta(\frac{z_{1}-z_{2}}{x})$ is defined in Section 2.

In order to make this work reasonably self-contained, we include
elementary definitions and notions needed in both theories.  The
reader can consult \cite{FLM} and \cite{FHL}, for example, for the
motivation and development of the theory of vertex operator algebras,
and \cite{Ha} for sheaves and \cite{Borel} for algebraic ${\cal
D}$-modules, whose theory was developed by Beilinson and Bernstein.

This paper is organized as follows: In Section 2, we recall some basic
notations and elementary tools and give the definitions of vertex
algebra and vertex algebra without vacuum. In Section 3, we recall
some basic concepts in the theory of ${\cal D}$-modules and give
examples which we shall need later. Beilinson-Drinfeld's notion of
chiral algebra over $X$ for a nonempty open subset $X\subset {\Bbb C}$
is given in Section 4. In Section 5, we define the notion of
vertex algebra without vacuum over $X$ and prove the equivalence theorem
stated above.

We would like to thank P. Deligne and especially A. Beilinson for
explaining the unpublished
work \cite{BD} to us.  We are also grateful to F. Knop,
whose Rutgers lecture notes on ${\cal D}$-modules and representation
theory were very helpful to us.  Y.-Z.~H. is supported in part by NSF
grant DMS-9596101 and by DIMACS, an NSF Science and Technology Center
funded under contract STC-88-09648, and J.~L by NSF grant DMS-9401851.

\section{Vertex algebras and vertex algebras without vacuum}

Following the treatment in \cite{FLM} and \cite{FHL}, we describe the
basic notations and elementary tools needed to formulate the notion of
vertex algebra.  We work over ${\Bbb C}$.
In this paper, the symbols $x, x_{0}, x_{1}, \dots$ are independent
commuting formal variables, and all expressions involving these
variables are to be understood as formal Laurent series. (Later we
shall also use the symbols $z, z_{1}, \dots,$ which will denote
complex numbers, not formal variables.)  We use the ``formal
$\delta$-function''
\[\delta(x)=\sum_{n\in {\Bbb Z}}x^{n},\]
which has the following simple and
fundamental property:
For any Laurent polynomial $f(x)\in {\Bbb C}[x, x^{-1}]$,
\[f(x)\delta(x)=f(1)\delta(x).\]
This property has many important variants. For example, for any
\[X(x_{1},
x_{2})\in (\mbox{End }W)[[x_{1}, x_{1}^{-1}, x_{2}, x_{2}^{-1}]]\]
(where $W$ is a vector space) such that
\[
\lim_{x_{1}\to x_{2}}X(x_{1}, x_{2})=X(x_{1},
x_{2})\mbar_{x_{1}=x_{2}}
\]
exists, we have
\begin{equation}\label{2.1}
X(x_{1}, x_{2})\delta\left(\frac{x_{1}}{x_{2}}\right)=X(x_{2}, x_{2})
\delta\left(\frac{x_{1}}{x_{2}}\right).
\end{equation}
The existence of this ``algebraic limit'' means that
for an arbitrary vector $w\in W$, the coefficient of each power of
$x_{2}$ in the formal expansion $X(x_{1}, x_{2})w\mbar_{x_{1}=x_{2}}$
is a finite sum.

We use the convention that negative powers of a
binomial are to be expanded in nonnegative powers of the second
summand. For example,
\[
x_{0}^{-1}\delta\left(\frac{x_{1}-x_{2}}{x_{0}}\right)=\sum_{n\in {\Bbb Z}}
\frac{(x_{1}-x_{2})^{n}}{x_{0}^{n+1}}=\sum_{m\in {\Bbb N},\; n\in {\Bbb Z}}
(-1)^{m}{{n}\choose {m}} x_{0}^{-n-1}x_{1}^{n-m}x_{2}^{m}.
\]
We have the following elementary identities:
\begin{eqnarray}
&{\dps x_{1}^{-1}\delta\left(\frac{x_{2}+x_{0}}{x_{1}}\right)
=x_{2}^{-1}\left(
\frac{x_{1}-x_{0}}{x_{2}}\right),}&\\
&{\dps x_{0}^{-1}\delta\left(\frac{x_{1}-x_{2}}{x_{0}}\right)-
x_{0}^{-1}\delta\left(\frac{x_{2}-x_{1}}{-x_{0}}\right)=
x_{2}^{-1}\delta\left(\frac{x_{1}-x_{0}}{x_{2}}\right).\label{s3}}&
\end{eqnarray}
Here and below, it is important to note
that the
relevant sums and products, etc., of formal series, are well defined. See
\cite{FLM} and \cite{FHL} for extensive discussions of the calculus
of formal $\delta$-functions.

The following version of the definition of vertex algebra, using
formal variables and the Jacobi identity of \cite{FLM} is equivalent
to Borcherds' original definition \cite{B}:

\begin{defi}
{\rm A {\it vertex algebra} is a
vector space $V$,
equipped with a linear map  $V\otimes V\to V[[x, x^{-1}]]$, or
equivalently,
\begin{eqnarray*}
V&\to&(\mbox{\rm End}\; V)[[x, x^{-1}]]\nonumber \\
v&\mapsto& Y(v, x)={\displaystyle \sum_{n\in{\Bbb Z}}}v_{n}x^{-n-1}
\;\;(\mbox{\rm where}\; v_{n}\in
\mbox{\rm End} \;V),
\end{eqnarray*}
$Y(v, x)$ denoting the {\it vertex operator associated with} $v$, and
equipped also with a distinguished homogeneous vector ${\bf 1}\in
V_{(0)}$ (the {\it vacuum}) and a linear map $D: V\to V$. The following
conditions are assumed for $u, v \in V$: the {\it lower truncation
condition} holds:
$$
u_{n}v=0\;\;\mbox{\rm for}\;n\; \mbox{\rm sufficiently large}
$$
(or equivalently, $Y(u, x)v\in V((x))$);
$$
Y({\bf 1}, x)=1\;\; (1\;\mbox{\rm on the right being the identity
operator});
$$
the {\it creation property} holds:
$$
Y(v, x){\bf 1} \in V[[x]]\;\;\mbox{\rm and}\;\;\lim_{x\rightarrow
0}Y(v, x){\bf 1}=v
$$
(that is, $Y(v, x){\bf 1}$ involves only nonnegative integral powers
of $x$ and the constant term is $v$); the {\it Jacobi identity} (the
main axiom) holds:
\begin{eqnarray*}
&x_{0}^{-1}\delta
\left({\displaystyle\frac{x_{1}-x_{2}}{x_{0}}}\right)Y(u, x_{1})Y(v,
x_{2})-x_{0}^{-1} \delta
\left({\displaystyle\frac{x_{2}-x_{1}}{-x_{0}}}\right)Y(v, x_{2})Y(u,
x_{1})&\nonumber \\ &=x_{2}^{-1} \delta
\left({\displaystyle\frac{x_{1}-x_{0}}{x_{2}}}\right)Y(Y(u, x_{0})v,
x_{2})&
\end{eqnarray*}
(note that when each expression in (2.8) is applied to any element of
$V$, the coefficient of each monomial in the formal variables is a
finite sum; on the right-hand side, the notation $Y(\cdot, x_{2})$ is
understood to be extended in the obvious way to $V[[x_{0},
x^{-1}_{0}]]$);
and
$$
{\displaystyle \frac{d}{dx}}Y(v,
x)=Y(Dv, x)
$$
(the {\it  $D$-derivative property}).}
\end{defi}

The vertex algebra just defined is denoted by $(V, Y, {\bf 1}, D)$
(or simply by $V$).   Homomorphisms of vertex
algebras are defined in the obvious way.

A consequence of the
definition above is the {\it skew-symmetry} \cite{B}:
\begin{equation}\label{skew}
Y(u, x)v=e^{xD}Y(v, -x)u
\end{equation}
for $u, v\in V$.  The proof uses the Jacobi identity, properties of
the $\delta$-function, the $D$-derivative property and the creation
property (see \cite{FHL}).

In the definition above, we have required that the vertex algebra has
a vacuum, but in Section 5 we shall see that the notion of chiral
algebra formulated using ${\cal D}$-modules does not give a vacuum.
Thus we need the following notion:

\begin{defi}\label{va-no-va}
{\rm A {\it vertex algebra without vacuum} is a vector space $V$, a vertex
operator map $Y: V\otimes V\to V[[x, x^{-1}]]$ and an operator $D$ on
$V$ satisfying all the axioms for a vertex algebra
 which do not involve the vacuum and in addition
the skew-symmetry (\ref{skew}).}
\end{defi}

We denote the vertex algebra without vacuum by $(V, Y, D)$ or
simply by $V$.

\renewcommand{\theequation}{\thesection.\arabic{equation}}
\renewcommand{\therema}{\thesection.\arabic{rema}}
\setcounter{equation}{0}
\setcounter{rema}{0}

\section{${\cal D}$-modules}

 We recall the elementary notions in the theories of sheaves (see
\cite{Ha}) and ${\cal D}$-modules (see \cite{Borel}).  We begin with
the definition of presheaf.  Following \cite{Borel}, we shall work
over nonsingular quasi-projective varieties over ${\Bbb C}$ of pure
dimension.  Below, by a variety we shall mean a variety of this type.

\begin{defi}
{\rm Let $X$ be a variety. A {\it presheaf ${\cal F}$ of abelian
groups on $X$}
consists of the  data

\begin{enumerate}

\item for every open subset $U\subset X$, an abelian group ${\cal
F}(U)$,

\item for every inclusion $V\subset U$ of open subsets of $X$, a
morphism of abelian groups $\rho_{UV}: {\cal F}(U)\to {\cal F}(V)$,

\end{enumerate}

\noindent satisfying the conditions

\begin{enumerate}

\item ${\cal F}(\emptyset)=0$,

\item $\rho_{UU}$ is the identity map {}from ${\cal F}(U)$ to itself,

\item if $W\subset V\subset U$ are open subsets of $X$, then
$\rho_{UW}=\rho_{VW}\circ \rho_{UV}$.

\end{enumerate}

\noindent We call the elements of ${\cal F}(U)$ {\it sections over
$U$} and the maps $\rho_{UV}$ the {\it restriction maps}.  If ${\cal
F}$ is a presheaf on $X$, and if $P$ is a point of $X$, we define the
{\it stalk} ${\cal F}_{P}$ of ${\cal F}$ at $P$ to be the direct limit
of the groups ${\cal F}(U)$ for all open sets $U$ containing $P$ via
the restriction maps $\rho$. An element of ${\cal F}_{P}$ is called a
{\it germ of sections of ${\cal F}$ at the point $P$}. }
\end{defi}

{\it Presheaves with values in any fixed category} can be defined by
replacing the references to abelian groups in the definition by the
analogous references for the given category. Most of the
considerations below hold for {\it presheaves of (commutative or
noncommutative) rings}, {\it presheaves of vector spaces} and {\it
presheaves of (commutative or noncommutative) algebras}, for example;
when we consider the zero object or subobjects or quotient objects, it
will be understood that we are always working in a category where
these exist.  We
 sometimes write $s|_{V}$ instead of $\rho_{UV}(s)$ for $s\in {\cal
F}(U)$. When it is necessary to distinguish the restriction maps of
different sheaves, we shall write the restriction maps for a sheaf
${\cal F}$ as $\rho^{\cal F}_{UV}$.

\begin{defi}

{\rm A presheaf on a variety $X$ is a {\it sheaf\/} if it satisfies
the following additional conditions:

\begin{enumerate}

\item If $U$ is an open subset of $X$, $\{V_{i}\}$ an open
covering of $U$ and $s\in {\cal F}(U)$ satisfying $s|_{V_{i}}=0$
for all $i$, then $s=0$.

\item If $U$ is an open subset of $X$, $\{V_{i}\}$ an open
covering of $U$ and we have elements $s_{i}\in {\cal F}(V_{i})$ for
each $i$, with the property that for each $i, j$, $s_{i}|_{V_{i}\cap
V_{j}}=s_{j}|_{V_{i}\cap V_{j}}$, then there is an element $s\in {\cal
F}(U)$ such that $s_{V_{i}}=s_{i}$ for each $i$.

\end{enumerate}}
\end{defi}

\begin{exam}\label{regfun}
{\rm Let $X$ be a variety.  We define the {\it sheaf ${\cal O}_{X}$ of
regular functions on $X$} as follows: For any open set $U$ of $X$, let
${\cal O}_X(U)$ be the ring of regular functions on $U$, which we
shall also write as ${\cal O}(U)$, and for any open subset $V\subset
U$ and $f\in {\cal F}(U)$, let $\rho_{UV}(f)$ be the restriction of
$f$ to $V$. It is clear that this gives a sheaf of rings (actually of
commutative associative algebras).}
\end{exam}

\begin{defi}
{\rm {\it Morphisms} and {\it isomorphisms of presheaves} are defined
in the obvious way. {\it Morphisms} and {\it isomorphisms of sheaves}
are defined to be morphisms and isomorphisms of the underlying
presheaves. {\it Presheaf kernels}, {\it presheaf cokernels} and {\it
presheaf images of morphisms of presheaves} are defined in the obvious
ways.}
\end{defi}

Note that presheaf kernels of morphisms of sheaves are sheaves
but in general presheaf cokernels and presheaf images of
morphisms of sheaves are not sheaves.

The following result associates a sheaf naturally to a presheaf in
terms of a universal property:

\begin{propo}
Given a presheaf ${\cal F}$, there is a sheaf ${\cal F}^{+}$ and a
morphism $\theta: {\cal F}\to {\cal F}^{+}$ with the property that for
any sheaf ${\cal G}$ and any morphism $\varphi: {\cal F}\to {\cal G}$,
there is a unique morphism $\psi: {\cal F}^{+}\to {\cal G}$ such that
$\varphi=\psi\circ \theta$. The pair $({\cal F}^{+}, \theta)$ is
unique up to unique isomorphism.
\end{propo}

See \cite{Ha}, for example, for a proof. The sheaf ${\cal F}^{+}$ is
called the {\it sheaf associated with the presheaf} ${\cal F}$.  It is
determined completely by the germs of sections of ${\cal F}$. Thus we
shall also call ${\cal F}^{+}$ the {\it sheaf of germs of sections of
${\cal F}$}.

\begin{defi}
{\rm A {\it subsheaf} of a sheaf ${\cal F}$ is a sheaf ${\cal F}'$
such that for every open set $U\subset X$, ${\cal F}'(U)$ is a
subgroup (or subobject in the category we are discussing) of ${\cal
F}(U)$, and the restriction maps of ${\cal F}'$ are induced {}from those
of ${\cal F}$.  Let $\varphi: {\cal F}\to {\cal G}$ be a morphism of
sheaves. The {\it kernel} of $\varphi$ is the presheaf kernel of
$\varphi$, which is a sheaf, and in fact a subsheaf. The {\it image}
of $\varphi$ is the sheaf associated to the presheaf image of
$\varphi$. The image of $\varphi$ can be identified with a subsheaf of
${\cal G}$ by means of the universal property of the sheaf associated
to a presheaf.  Let ${\cal F}'$ be a subsheaf of a sheaf ${\cal
F}$. The {\it quotient sheaf} ${\cal F}/{\cal F}'$ is the sheaf
associated to the presheaf given by $U\mapsto {\cal F}(U)/{\cal F}'(U)$.
The {\it cokernel} of $\varphi$ is the sheaf associated to the
presheaf cokernel of $\varphi$. }
\end{defi}

\begin{defi}\label{3.7}

{\rm Let $X$ and $Y$ be varieties and $f: X \to Y$ a morphism. For any
sheaf ${\cal F}$ on $X$, The {\it direct image} $f_{*}({\cal F})$ is
the sheaf on $Y$ given by $V\mapsto (f_{*}({\cal F}))(V)={\cal
F}(f^{-1}(V))$ for any open set $V\subset Y$. For any sheaf ${\cal G}$
on $Y$, the {\it inverse image} $f^{-1}({\cal G})$ is the sheaf associated to
the presheaf given by $U\mapsto \lim_{V \supset f(U)}{\cal G}(V)$, where
$U$ is any open set of $X$ and the limit is taken over all open sets
$V$ of $Y$ containing $f(U)$. }
\end{defi}

For any open subset $V\subset Y$ and $s\in {\cal
F}(f^{-1}(V))$, we denote $s$ by $f_{*}(s)$ when it is viewed as
an element of $(f_{*}({\cal F}))(V)$, and we shall call
$s$ the {\it preimage of $f_{*}(s)$}. For any open subset 
$U\subset X$ and any
open subset $V\subset Y$ such that $V\supset f(U)$, an element
$s\in {\cal G}(V)$ determines an element of $(f^{-1}({\cal G}))(U)$.
We denote this element by $f^{-1}(s)$.

We now define the notions of left and right ${\cal O}$-module for any
sheaf ${\cal O}$ of rings.  Modules will always be left modules unless
``right module'' is specified.  The concepts below have obvious
analogues for right modules.

\begin{defi}\label{omodule}
{\rm Let ${\cal O}$ be a sheaf of rings on $X$.
A sheaf ${\cal F}$ on $X$ is an {${\cal O}$-module} if for
any open set $U$, ${\cal F}(U)$ is an ${\cal O}(U)$-module,
and for
any open sets $V\subset U$ and $f\in {\cal O}(U)$, $s\in {\cal F}(U)$,
we have $\rho^{\cal F}_{UV}(fs)=\rho^{{\cal O}}_{UV}(f)\rho^{\cal
F}_{UV}(s)$.}
\end{defi}

{\it Morphisms} and {\it isomorphisms of ${\cal O}$-modules} are
defined in the obvious ways. Note that kernels, images and cokernels
of morphisms of ${\cal O}$-modules are again ${\cal O}$-modules. If
${\cal F}'$ is a subsheaf of ${\cal O}$-modules of an ${\cal
O}$-module ${\cal F}$, then the quotient sheaf ${\cal F}/{\cal F}'$ is
also an ${\cal O}$-module.

Now we restrict our attention to the sheaf ${\cal O}_X$ (recall
Example \ref{regfun}), for which left and right modules are the same.

\begin{defi}
{\rm An ${\cal O}_X$-module ${\cal F}$ is said to be {\it
quasi-coherent} if for any open affine subset $U$ of $X$ and $f\in
{\cal O}_X(U)$, the following conditions hold:

\begin{enumerate}

\item Let $V$ be the open set $\{x\in U\;|\: f(x)\ne 0\}$. For any
$s\in {\cal F}(V)$, there exists $n\in {\Bbb N}$ and $\bar{s}\in {\cal
F}(U)$ such that $\rho_{UV}(\bar{s})=f^{n}s$.

\item For any $s\in {\cal F}(U)$ such that $\rho_{UV}(s)=0$, there
exists $n\in {\Bbb N}$ such that $f^{n}s=0$.

\end{enumerate}}

\end{defi}

\begin{defi}
{\rm Let ${\cal F}$ and ${\cal G}$ be ${\cal O}_X$-modules. The
{\it tensor product} ${\cal F}\otimes_{{\cal O}_X}{\cal G}$ (or simply
${\cal F}\otimes {\cal G}$) is the sheaf associated to the presheaf
given by $U\mapsto {\cal F}(U)\otimes_{{\cal O}_X(U)}{\cal G}(U)$.}
\end{defi}

\begin{defi}\label{3.11}
{\rm Let $X$ and $Y$ be varieties, $f: X\to Y$ a morphism.  This
induces a natural morphism $f^{\#}: {\cal O}_Y \to f_{*}({\cal O}_X)$ of
sheaves. If ${\cal F}$ is an ${\cal O}_X$-module, then $f_{*}({\cal F})$
is an $f_{*}({\cal O}_X)$-module, and thus an ${\cal O}_Y$-module via
$f^{\#}$.  This ${\cal O}_Y$-module $f_{*}({\cal F})$ is called the {\it
direct image of ${\cal F}$ by $f$}.  If ${\cal G}$ is a ${\cal
O}_Y$-module, then $f^{-1}({\cal G})$ is an $f^{-1}({\cal
O}_Y)$-module. Since $f^{\#}$ induces a natural morphism {}from
$f^{-1}({\cal O}_Y)$ to ${\cal O}_X$, we can form the tensor product
$f^{-1}({\cal G})\otimes_{f^{-1}({\cal O}_Y)}{\cal O}_X$ and it is an
${\cal O}_X$-module.  This ${\cal O}_X$-module is called the {\it
inverse image of ${\cal G}$ by $f$} and is denoted $f^{*}({\cal G})$.}
\end{defi}

For any open subset $U\subset X$, any
open subset $V\subset Y$ such that $V\supset f(U)$, and any element
$s\in {\cal G}(V)$, we denote
$f^{-1}(s)\otimes_{(f^{-1}({\cal O}_Y))(U)} 1\in (f^{*}({\cal G}))(U)$
by $f^{*}(s)$.

\begin{exam}
{\rm Let $X$ be a variety.  Here we define the sheaf of
(noncommutative) algebras ${\cal D}_X$---the {\it sheaf of germs of
algebraic differential operators on $X$} (see \cite{Borel}).  If $X$
is an affine variety, we define ${\cal D}_{X}(X)$ to be the algebra
${\cal D}(X)$ of operators on ${\cal O}(X)={\cal O}_{X}(X)$ generated
by the elements of ${\cal O}(X)$ (acting by multiplication) and the
derivations of ${\cal O}(X)$. This is the {\it ring of (algebraic)
differential operators on $X$}. For any open affine subset $U$ of $X$,
we have ${\cal D}(U)= {\cal O}(U)\otimes_{{\cal O}(X)}{\cal D}(X)$,
and there exists a unique quasi-coherent ${\cal O}_{X}$-module ${\cal
D}_X$ on $X$ such that ${\cal D}(U)$ is the ${\cal O}(U)$-module of
sections over an open affine subset $U\subset X$. In general, for any
variety $X$ (of the type we are considering), there exists a unique
sheaf ${\cal D}_{X}$ of algebras on $X$ whose restriction to every
open affine subset $U\subset X$ is ${\cal D}_{U}$. The algebra ${\cal
D}_{X}(X)$ consists of the {\it algebraic differential operators on
$X$}. The sheaf ${\cal D}_{X}$ is quasi-coherent as an ${\cal
O}_{X}$-module and its restriction to an open subset $U\subset X$ is
${\cal D}_{U}$.}
\end{exam}

We shall use the following notion of (algebraic) ${\cal D}$-module:

\begin{defi}
{\rm Let $X$ be a variety. A {\it ${\cal D}_{X}$-module}  is a
${\cal D}_{X}$-module as defined above which is also
quasi-coherent as an ${\cal O}_{X}$-module. Similarly for a {\it right
${\cal D}_X$-module}.}
\end{defi}

Note that a (left or right) ${\cal D}_{X}$-module is a
 ${\cal D}_{X}$-module in the sense of Definition \ref{omodule}
satisfying an extra condition.
 {\it
Morphisms} and {\it isomorphisms of
${\cal D}_{X}$-modules} are defined in
the obvious ways. It is clear that the set of all morphisms
{}from a
${\cal D}_{X}$-module to another one has a natural
abelian group structure.
For any (left) ${\cal D}_{X}$-modules
${\cal F}_{1}, {\cal F}_{2}$ and right ${\cal D}_{X}$-modules
${\cal G}_{1}, {\cal G}_{2}$, ${\cal F}_{1}
\otimes _{{\cal O}_{X}}{\cal F}_{2}$ and
$\hom_{{\cal O}_{X}}({\cal G}_{1}, {\cal G}_{2})$
are naturally (left)  ${\cal D}_{X}$-modules, and
${\cal F}_{1}
\otimes _{{\cal O}_{X}}{\cal G}_{1}$ and
$\hom_{{\cal O}_{X}}({\cal F}_{1}, {\cal G}_{2})$
are naturally right ${\cal D}_{X}$-modules.
In particular,
we have a tensor product operation in the
category of (left) ${\cal D}_{X}$-modules.

\begin{exam}
{\rm Let $\omega_{X}$ be the sheaf of germs of differential forms of
top degree on $X$.  Let $U$ be an open affine subset of $X$. A
derivation $\xi\in {\cal D}_{X}(U)$ acts on $\omega_{X}(U)$ {}from the
left by the Lie derivative $L_{\xi}$, but this does not give a (left)
${\cal D}_{X}$-module structure to $\omega_{X}$.  Rather, the action
$-L_{\xi}$ of $\xi$ on $\omega_{X}(U)$ gives a right ${\cal
D}_{X}$-module structure to $\omega_{X}$. Thus for any (left) ${\cal
D}_{X}$-module ${\cal F}$ and right ${\cal D}_{X}$-module ${\cal G}$,
${\cal F}\otimes_{{\cal O}_{X}}\omega_{X}$ is a right ${\cal
D}_{X}$-module and $\hom_{{\cal O}_{X}}(\omega_{X},{\cal G})$ is a
(left) ${\cal D}_{X}$-module.  (See \cite{Borel}.)}
\end{exam}

We have the following notions of inverse image and direct image of a
${\cal D}_{X}$-module:

\begin{defi}
{\rm Let $X$ and $Y$ be varieties and $f: X\to Y$ a morphism.  For any
${\cal D}_{Y}$-module ${\cal G}$, we have the inverse image ${\cal
O}_{X}$-module $f^{*}({\cal G})$, which is quasi-coherent.  There is a
natural ${\cal D}_{X}$-module structure on $f^{*}({\cal G})$ (see
\cite{Borel}), which we call the {\it inverse image ${\cal
D}_{X}$-module} and write as $f^\circ({\cal G})$.  Now, to define the
direct image, we begin with the ${\cal O}_{X}$-module $f^{*}({\cal
D}_{Y})$, which has a natural (left) ${\cal D}_{X}$-module structure,
as above.  On the other hand, the right multiplication on ${\cal
D}_{Y}$ carries over to a right $f^{-1}({\cal D}_{Y})$-module structure
on $f^{*}({\cal D}_{Y})$.  We denote $f^{*}({\cal D}_{Y})$ by ${\cal
D}_{X\to Y}$ when it is equipped with this ${\cal D}_{X}\times
f^{-1}({\cal D}_{Y})$-module structure. Let ${\cal D}_{Y\leftarrow
X}={\cal D}_{X\to Y} \otimes_{{\cal O}_{X}}\hom_{{\cal
O}_{X}}(f^{*}(\omega_{Y}), \omega_{X})$.  Then ${\cal D}_{Y\leftarrow
X}$ is a (left) $f^{-1}({\cal D}_{Y})$-module and a right ${\cal
D}_{X}$-module. Let ${\cal F}$ be a (left) ${\cal D}_{X}$-module. Then
${\cal D}_{Y\leftarrow X} \otimes_{{\cal D}_{X}}{\cal F}$ (defined in
the obvious way) is a (left) $f^{-1}({\cal D}_{Y})$-module.  We also
have a canonical morphism $\tilde{f}: {\cal D}_{Y}\to
f_{*}(f^{-1}({\cal D}_{Y}))$ of sheaves of rings (recall Definition
\ref{3.7}).  We define the {\it direct image ${\cal D}_{Y}$-module}
$f_{\circ}({\cal F})$ to be $f_*({\cal D}_{Y\leftarrow X}
\otimes_{{\cal D}_{X}}{\cal F})$, viewed as a ${\cal D}_{Y}$-module
via $\tilde{f}$ (cf. the definition of direct image in Definition
\ref{3.11}).}
\end{defi}

\begin{exam}\label{delta*}
{\rm Let $X$ be a nonempty open subset of ${\Bbb C}$, $Y=X\times X$,
$f=\Delta: X\to Y=X\times X$ the diagonal map, defined by
$\Delta(z)=(z, z)$, and ${\cal F}$ a ${\cal D}_{X}$-module. Note that
$X$ is necessarily affine, and that for a local theory over curves,
such varieties $X$ are enough.  In this case, $\omega_{X}(X)$ and
$\omega_{X\times X}(X\times X)$ are generated by $dz$ and
$dz_{1}\wedge dz_{2}$, respectively, and
$$(\hom_{{\cal O}_{X}}(\Delta^{*}(\omega_{X\times X}), \omega_{X}))(X)=
\hom_{{\cal O}_{X}(X)}(\Delta^{*}(\omega_{X\times X}(X\times X)),
\omega_{X}(X))$$ is
generated by $\phi\in
\hom_{{\cal O}_{X}(X)}((\Delta^{*}(\omega_{X\times X}(X\times X)),
\omega_{X}(X))$
defined by
$$\phi(\Delta^{*}(dz_{1}\wedge dz_{2}))=
dz.$$
Thus $(\Delta_{\circ}({\cal F}))(X\times X)$
is generated by elements of the form
$$\Delta_{*}((\Delta^{*}(1)
\otimes_{{\cal O}_{X}(X)}\phi)
\otimes_{{\cal D}_{X}(X)}A)$$
for $A\in {\cal F}(X)$.
The element of $\frac{\p}{\p z_{1}}+\frac{\p}{\p z_{2}}$ of
${\cal D}_{X\times X}(X\times X)$ is in fact in
$(\Delta_{*}(\Delta^{-1}({\cal D}_{X\times X})))(X\times X)$
and {}from the definition, we have
\begin{eqnarray*}
\lefteqn{\left(\frac{\p}{\p z_{1}}+\frac{\p}{\p z_{2}}\right)
\Delta_{*}((\Delta^{*}(\xi)
\otimes_{{\cal O}_{X}(X)}\phi)
\otimes_{{\cal D}_{X}(X)}A)}\nno\\
&&=\Delta_{*}\left(\left(\Delta^{*}
\left(\xi\left(\frac{\p}{\p z_{1}}+\frac{\p}{\p z_{2}}\right)\right)
\otimes_{{\cal O}_{X}(X)}\phi\right)
\otimes_{{\cal D}_{X}(X)}A\right)\nno\\
&&=\Delta_{*}\left(\left(\Delta^{*}(\xi)
\otimes_{{\cal O}_{X}(X)}\phi\right)
\otimes_{{\cal D}_{X}(X)}\frac{\p}{\p z}A\right)
\end{eqnarray*}
for any $\xi\in {\cal D}_{X\times X}(X\times X)$ commuting
with $\frac{\p}{\p z_{1}}+\frac{\p}{\p z_{2}}$ and $A\in {\cal F}(X)$,
so that the elements
$$\frac{\p^{n}}{\p z_{1}^{n}}\Delta_{*}\left(\left(\Delta^{*}
\left(1\right)
\otimes_{{\cal O}_{X}(X)}\phi\right)
\otimes_{{\cal D}_{X}(X)}A\right),\;\;n\in \Bbb{Z},$$
span $(\Delta_{\circ}({\cal F}))(X\times X)$.
Note that since $X$ is an open subset of ${\Bbb C}$,
$(\Delta_{*}({\cal D}_{X}))(X\times X)$ can be embedded
into ${\cal D}_{X\times X}(X\times X)$ such that
$\Delta_{*}(\frac{\p}{\p z})$ is mapped to
$\frac{\p}{\p z_{1}}+\frac{\p}{\p z_{2}}$. Thus we have a
${\cal D}_{X\times X}(X\times X)$-module
$${\cal D}_{X\times X}(X\times X)
\otimes_{(\Delta_{*}({\cal D}_{X}))(X\times X)}
(\Delta_{*}({\cal F}))(X\times X),$$
and if $\{A_{\alpha}\}$ is a basis of ${\cal F}(X)$, this
${\cal D}_{X\times X}(X\times X)$-module has as a basis
$\frac{\p^{n}}{\p z_{1}^{n}}
\otimes_{(\Delta_{*}({\cal D}_{X}))(X\times X)} A_{\alpha}$ for
$n\in {\Bbb N}$.
It is easy to see that
$$\frac{\p^{n}}{\p z_{1}^{n}}\Delta_{*}\left(\left(\Delta^{*}
\left(1\right)
\otimes_{{\cal O}_{X}(X)}\phi\right)
\otimes_{{\cal D}_{X}(X)}A\right)\mapsto
\frac{\p^{n}}{\p z_{1}^{n}}
\otimes_{(\Delta_{*}({\cal D}_{X}))(X\times X)} A$$
gives an isomorphism {}from
 $(\Delta_{\circ}({\cal F}))(X\times X)$  to
$${\cal D}_{X\times X}(X\times X)
\otimes_{(\Delta_{*}({\cal D}_{X}))(X\times X)}
(\Delta_{*}({\cal F}))(X\times X).$$
Thus $(\Delta_{\circ}({\cal F}))(X\times X)$ has as a basis
the elements of the form
$$\frac{\p^{n}}{\p z_{1}^{n}}\Delta_{*}\left(\left(\Delta^{*}
\left(1\right)
\otimes_{{\cal O}_{X}(X)}\phi\right)
\otimes_{{\cal D}_{X}(X)}A_{\alpha}\right)$$
for $n\in {\Bbb N}$, and we can
identify the ${\cal D}_{X\times X}$-module $\Delta_{\circ}({\cal F})$ with
${\cal D}_{X\times X}\otimes_{\Delta_{*}({\cal D}_{X})}
\Delta_{*}({\cal F})$.
Now let us take ${\cal F}={\cal O}_{X}$. Then we have
the element
\begin{equation}\label{d-delta}
\Delta_{*}((\Delta^{*}(1)
\otimes_{{\cal O}_{X}(X)}\phi)
\otimes_{{\cal D}_{X}(X)}1).
\end{equation}
For any $f(z_{1}, z_{2})\in {\cal O}_{X\times X}(X\times X)$,
we have
\begin{eqnarray}\label{d-delta-limit}
\lefteqn{f(z_{1}, z_{2})\Delta_{*}((\Delta^{*}(1)
\otimes_{{\cal O}_{X}(X)}\phi)
\otimes_{{\cal D}_{X}(X)}1)}\nno\\
&&=\Delta_{*}(f(z, z)((\Delta^{*}(1)
\otimes_{{\cal O}_{X}(X)}\phi)
\otimes_{{\cal D}_{X}(X)}1))\nno\\
&&=f(z_{1}, z_{1})\Delta_{*}((\Delta^{*}(1)
\otimes_{{\cal O}_{X}(X)}\phi)
\otimes_{{\cal D}_{X}(X)}1).
\end{eqnarray}
Similarly, the left-hand side of (\ref{d-delta-limit}) is also equal to
$$f(z_{2}, z_{2})\Delta_{*}((\Delta^{*}(1)
\otimes_{{\cal O}_{X}(X)}\phi)
\otimes_{{\cal D}_{X}(X)}1).$$
So we see that the element (\ref{d-delta}) has the property similar
to  the property (\ref{2.1}) of the $\delta$-function
$x^{-1}_{2}\delta\left(\frac{x_{1}}{x_{2}}\right)$.
In fact it is easy to show that the derivatives
$$\frac{\p^{m}}{\p z_{1}^{m}}\Delta_{*}((\Delta^{*}(1)
\otimes_{{\cal O}_{X}(X)}\phi)
\otimes_{{\cal D}_{X}}1),$$
$m\in {\Bbb N}$, and
$$\frac{\p^{n}}{\p z_{2}^{n}}\Delta_{*}((\Delta^{*}(1)
\otimes_{{\cal O}_{X}(X)}\phi)
\otimes_{{\cal D}_{X}}1),$$
$n\in {\Bbb N}$, also have the corresponding properties of the derivatives
of $x^{-1}_{2}\delta\left(\frac{x_{1}}{x_{2}}\right)$ (see
\cite{FLM}, Proposition 8.2.2). Thus we can identify
(\ref{d-delta}) with $x^{-1}_{2}\delta\left(\frac{x_{1}}{x_{2}}\right)$.
In particular, we obtain a ${\cal D}_{X\times X}$-module
generated by $x^{-1}_{2}\delta\left(\frac{x_{1}}{x_{2}}\right)$ which
is isomorphic to $\Delta_{\circ}({\cal O}_{X})$.}
\end{exam}

\begin{rema}\label{brokens3}
{\rm In the example above, we have given a ${\cal D}$-module-theoretic
interpretation of the $\delta$-function
$x^{-1}_{2}\delta\left(\frac{x_{1}}{x_{2}}\right)$ and its derivatives.
But in the formal variable approach to vertex algebras,
it is most natural to use
 $x^{-1}_{2}\delta\left(\frac{x_{1}-x_{0}}{x_{2}}\right)$ (see above),
which is in fact
a formal infinite linear combination  of derivatives of
$x^{-1}_{2}\delta\left(\frac{x_{1}}{x_{2}}\right)$, by a formal Taylor's
theorem (see e.g. \cite{FLM}, Propositions 8.2.2 and 8.3.1).
This $\delta$-function and the other two in the  Jacobi
identity for a vertex algebra  apparently
 do not have direct interpretations in
terms of ${\cal D}$-modules. The $S_{3}$-symmetry
of (\ref{s3}) and of the Jacobi identity for vertex algebras
(see \cite{FLM}, \cite{FHL})
involves all the three formal variables $x_{0}, x_{1}, x_{2}$ on equal
footing, and
we note that the ${\cal D}$-module
approach does not naturally have three variables
on equal footing.}
\end{rema}

We also need the concept of exterior tensor product:

\begin{defi}
{\rm Let $X$ and $Y$ be varieties and ${\cal F}$ and ${\cal G}$ ${\cal
D}_{X}$- and ${\cal D}_{Y}$-modules, respectively. Let $({\cal
F}\boxtimes {\cal G})(U\times V) ={\cal F}(U)\otimes {\cal G}(V)$ for
all open affine subsets $U\subset X$ and $V\subset Y$.  Then there is
a unique ${\cal D}_{X\times Y}$-module ${\cal F}\boxtimes {\cal G}$
such that $({\cal F}\boxtimes {\cal G})(U_{1}\times U_{2})$ is the
${\cal O}_{X \times Y}(U \times V)$-module of sections on $U\times V$
for any open affine subsets $U\subset X$ and $V\subset Y$.  We call
this ${\cal D}_{X\times Y}$-module the {\it exterior tensor product of
${\cal F}$ and ${\cal G}$}.  {}From the definition, we see that ${\cal
F}\boxtimes {\cal G}$ is determined by sections on open affine subsets
of the form $U\times V$.}
\end{defi}

\renewcommand{\theequation}{\thesection.\arabic{equation}}
\renewcommand{\therema}{\thesection.\arabic{rema}}
\setcounter{equation}{0} 
\setcounter{rema}{0}

\section{Beilinson-Drinfeld's chiral algebra}

We now give Beilinson-Drinfeld's definition of chiral algebra.  For
simplicity and for the purpose of comparing it with the definition of
vertex algebra without vacuum, we only give the definition of chiral
algebra over a nonempty open subset $X\subset {\Bbb C}$.  But since
 the definition
is sheaf-theoretic and thus is local in nature, this definition
carries over naturally to the general case.  The axioms in the
definition are natural in the language of quasi-tensor categories
in Beilinson-Drinfeld's formulation.

Let $X\subset {\Bbb C}$ be a nonempty open subset.
Let $F^{2}$ be the complement of the diagonal in
$X\times X$. We denote the embedding map {}from
$F^{2}$ to $X\times X$ by $j$ and the
diagonal map {}from $X$ to $X\times X$ by
$\Delta$. Let
$$F^{3}=\{(z_{1}, z_{2}, z_{3})\in X\times
X\times X\;|\; z_{k}\ne z_{l} \;\mbox{\rm for}\;
k\ne l, k, l=1, 2, 3\}.$$
We denote the embedding map {}from $F^{3}$ to $X\times
X\times X$ by $j_{3}$ and the diagonal map {}from
$X$ to $X\times
X\times X$ by $\Delta_{3}$.

Let ${\cal A}$ be a  ${\cal D}_{X}$-module.
The exterior tensor product ${\cal D}_{X\times X}$-module
${\cal A}\boxtimes {\cal A}$ is by definition determined by
\begin{eqnarray*}
\lefteqn{U_{1}\times U_{2}\subset X\times X}\nno\\
&&\mapsto \{\sum_{i=1}^{n}f_{i}
(A_{i}\otimes
B_{i})
\;|\;A_{i}\in {\cal A}(U_{1}), B_{i}\in {\cal A}(U_{2}),\nno\\
&&\hspace{10em}f_{i}\in {\cal O}(U_{1}\times U_{2})
\;\mbox{\rm for}\;i=1, \dots, n\}
\end{eqnarray*}
for all open subsets $U_{1}, U_{2}\subset X$.
The inverse image ${\cal O}_{F^{2}}$-module
$j^{*}({\cal A}\boxtimes {\cal A})$ is
 the ${\cal O}_{F^{2}}$-module obtained by restricting the
sections
of ${\cal A}\boxtimes {\cal A}$ over open subsets of $X\times X$
 to open subsets of $F^{2}$, and the direct image ${\cal O}_{X\times
X}$-module
$j_{*}j^{*}({\cal A}\boxtimes {\cal A})$ is
the ${\cal O}_{X\times
X}$-module determined by
\begin{eqnarray*}
\lefteqn{U_{1}\times U_{2}\subset X\times X}\nno\\
&&\mapsto \{\sum_{i=1}^{n}f_{i}
(A_{i}\otimes
B_{i})
\;|\;A_{i}\in {\cal A}(U_{1}), B_{i}\in {\cal A}(U_{2}), \nno\\
&&\hspace{10em}
f_{i}\in {\cal O}((U_{1}\times U_{2})\cap F^{2})
\;\mbox{\rm for}\;i=1, \dots, n\}
\end{eqnarray*}

for all open subsets $U_{1}, U_{2}\subset X$. It has a natural
${\cal D}_{X\times X}$-module structure.
As an ${\cal O}_{X\times
X}(U_{1}\times U_{2})$-module,
$j_{*}j^{*}({\cal A}\boxtimes {\cal A})(U_{1}\times U_{2})$ is
generated by the elements of the form
$$(z_{1}-z_{2})^{m}(A_{1}\otimes A_{2})$$
for $m<0$,
$A_{1}\in {\cal A}(U_{1})$ and $A_{2}\in {\cal A}(U_{2})$.
Similarly, if $U=U_{1}\times U_{2}\times U_{3}$ where $U_{1}, U_{2},
U_{3}$
are open subsets of $X$, then
$((j_{3})_{*}j_{3}^{*}({\cal A}\boxtimes {\cal A}\boxtimes {\cal A}))(U)$
as an ${\cal O}_{X\times
X\times X}(U)$-module is generated by the
elements of the form
$$(z_{1}-z_{2})^{m_{1}}(z_{2}-z_{3})^{m_{2}}(z_{1}-z_{3})^{m_{3}}
(A_{1}\otimes A_{2}\otimes A_{3})$$
for $m_{i}<0$, $A_{i}\in {\cal A}(U_{i})$, $i=1, 2, 3$. This
has a natural ${\cal D}_{X\times
X\times X}(U)$-module structure, and so we have the
${\cal D}_{X\times X\times X}$-module
$(j_{3})_{*}j_{3}^{*}({\cal A}\boxtimes {\cal A}\boxtimes {\cal A})$.

In the definition of Lie algebra, we need compositions
$[[\cdot, \cdot], \cdot]$ and $[\cdot, [\cdot, \cdot]]$  of the Lie
bracket $[\cdot, \cdot]$
to formulate the Jacobi identity. Here we need analogues of
these compositions.

For any  morphism
$\mu: j_{*}j^{*}({\cal A}\boxtimes {\cal A})\to \Delta_{\circ}{\cal A}$
of ${\cal D}_{X\times X}$-modules,
we define
a natural morphism $$\mu(\mu(\cdot, \cdot), \cdot):
(j_{3})_{*}j_{3}^{*}({\cal A}\boxtimes {\cal A}\boxtimes {\cal A})
\to (\Delta_{3})_{\circ}({\cal A})$$
of ${\cal D}_{X\times X\times X}$-modules as follows:
Let $U$ be an open subset of $X\times
X\times X$ of the form $U_{1}\times U_{2}\times U_{3}$
where $U_{1}, U_{2}$ and $U_{3}$ are open subsets of $X$.
We identify $\Delta_{\circ}({\cal A})$ with
${\cal D}_{X\times X}\otimes_{(\Delta_{*}({\cal D}_{X}))}
\Delta_{*}({\cal A})$ (see Example \ref{delta*}).
In particular,
$(\Delta_{\circ}({\cal A}))(U_{1}\times U_{2})$ is spanned by the
elements of the form
$$\frac{\p^{n}}{\p z_{1}^{n}}\otimes_{(\Delta_{*}({\cal D}_{X}))(U)}
\Delta_{*}(A)$$
for $n\in \Bbb{N}$ and
$A\in {\cal A}(\Delta^{-1}(U_{1}\times U_{2}))$. Thus
we see that for any $n_{1}\in {\Bbb Z}$ there exist
$p_{n_{1}}\in {\Bbb N}$ and
$$B^{n_{1}}_{k}\in {\cal A}(\Delta^{-1}(U_{1}\times U_{2})) ={\cal
A}(U_{1}\cap U_{2}), \;\;k=0, \dots, p_{n_{1}},$$
such that
\begin{equation}\label{expl-exp}
\mu((z_{1}-z_{2})^{n_{1}}(A_{1}\otimes A_{2}))
=\sum_{k=0}^{p_{n_{1}}}\frac{\p^{k}}{\p z^{k}_{1}}\otimes_{(\Delta_{*}({\cal
D}_{X}))(U_{1}\times U_{2})} \Delta_{*}(B^{n_{1}}_{k});
\end{equation}
the nonzero elements $B_{k}^{n_{1}}$ are uniquely determined.
Note that when $n_{1}\ge 0$, $(z_{1}-z_{2})^{n_{1}}$
is regular, so that when $n_{1}>p_{0}$, we have
\begin{eqnarray*}
\lefteqn{\mu((z_{1}-z_{2})^{n_{1}}(A_{1}\otimes A_{2}))=}\nno\\
&&=(z_{1}-z_{2})^{n_{1}}\mu(A_{1}\otimes A_{2})\nno\\
&&=(z_{1}-z_{2})^{n_{1}}
\sum_{k=0}^{p_{0}}\frac{\p^{k}}{\p z^{k}_{1}}\otimes_{(\Delta_{*}({\cal
D}_{X}))(U_{1}\times U_{2})} \Delta_{*}(B^{0}_{k})\nno\\
&&=0.
\end{eqnarray*}
 Similarly, for any $n_{2}\in {\Bbb Z}$ we have, for the elements
$B_{k}^{n_{1}}$,
\begin{equation}\label{4.0.1}
\mu((z_{2}-z_{3})^{n_{2}}(B^{n_{1}}_{k}\otimes A_{3}))
=\sum_{l=0}^{q_{n_{1}, n_{2}}}\frac{\p^{l}}{\p z^{l}_{2}}
\otimes_{(\Delta_{*}({\cal
D}_{X}))(U_{2}\times U_{3})} \Delta_{*}(C^{n_{1}n_{2}}_{kl})
\end{equation}
where $q_{n_{1}, n_{2}}\in {\Bbb N}$ and
$C^{n_{1}n_{2}}_{kl}\in {\cal A}(\Delta^{-1}((U_{1}\cap U_{2})
\times U_{3})
={\cal A}(U_{1}\cap U_{2}\cap U_{3})$,
$k=0, \dots, p_{n_{1}}$, $l=0, \dots, q_{n_{1}, n_{2}}$.
Since $\frac{\p}{\p z_{2}}+\frac{\p}{\p z_{3}}\in
(\Delta_{*}({\cal D}_{X}))(U_{2}\times U_{3})$ and its preimage in
${\cal D}_{X}(U_{2}\cap U_{3})$ is $\frac{\p}{\p z}$, the right-hand side of
(\ref{4.0.1}) is equal to
\begin{eqnarray}\label{4.0.2}
\lefteqn{\sum_{l=0}^{q_{n_{1}, n_{2}}}\left(\left(\frac{\p}{\p z_{2}}
+\frac{\p}{\p z_{3}}\right)-\frac{\p}{\p z_{3}}\right)^{l}
\otimes_{(\Delta_{*}({\cal
D}_{X}))(U_{2}\times U_{3})} \Delta_{*}(C^{n_{1}n_{2}}_{kl})}\nno\\
&&=\sum_{l=0}^{q_{n_{1}, n_{2}}}\sum_{j=0}^{l}
(-1)^{j}{l\choose j}\frac{\p^{j}}{\p z^{j}_{3}}\left(\frac{\p}{\p z_{2}}
+\frac{\p}{\p z_{3}}\right)^{l-j}
\otimes_{(\Delta_{*}({\cal
D}_{X}))(U_{2}\times U_{3})} \Delta_{*}(C^{n_{1}n_{2}}_{kl})\nno\\
&&=\sum_{l=0}^{q_{n_{1}, n_{2}}}\sum_{j=0}^{l}
(-1)^{j}{l\choose j}\frac{\p^{j}}{\p z^{j}_{3}}
\otimes_{(\Delta_{*}({\cal
D}_{X}))(U_{2}\times U_{3})}
\Delta_{*}\left(\frac{\p^{l-j}}{\p z^{l-j}}C^{n_{1}n_{2}}_{kl}\right).
\end{eqnarray}
Let $C^{n_{1}n_{2}}_{kl}=0$ for
$k>p_{n_{1}}$ or $l>q_{n_{1}, n_{2}}$.
As in Example \ref{delta*},
$(\Delta_{3})_{\circ}({\cal A})$ is canonically isomorphic to
${\cal D}_{X\times X\times X}\otimes_{(\Delta_{3})_{*}({\cal D}_{X})}
{\cal A}$, and we shall identify $(\Delta_{3})_{\circ}({\cal A})$ with
${\cal D}_{X\times X\times X}\otimes_{(\Delta_{3})_{*}({\cal D}_{X})}
{\cal A}$.
We define
\begin{eqnarray}\label{4.0.3}
 \lefteqn{(\mu(\mu(\cdot, \cdot), \cdot))
((z_{1}-z_{2})^{m_{1}}(z_{2}-z_{3})^{m_{2}}(z_{1}-z_{3})^{m_{3}}
(A_{1}\otimes A_{2}\otimes A_{3}))}\nno\\
&& =\sum_{i\in {\Bbb N}}
\sum_{k, l\in {\Bbb N}}\sum_{j=0}^{l}
{m_{3}\choose i}(-1)^{j}{l\choose j}\frac{\p^{k}}{\p z^{k}_{1}}
\frac{\p^{j}}{\p z^{j}_{3}}\nno\\
&&\quad\quad \otimes_{((\Delta_{3})_{*}
({\cal D}_{X}))(U)}(\Delta_{3})_{*}
\left(\frac{\p^{l-j}}{\p z^{l-j}}C^{m_{1}+i,
m_{2}+m_{3}-i}_{kl}\right),
\end{eqnarray}
using the expansion of
$$(z_{1}-z_{3})^{m_{3}}=((z_{2}-z_{3})+(z_{1}-z_{2}))^{m_{3}}$$
in nonnegative powers of $z_{1}-z_{2}$.
Note that the right-hand side of (\ref{4.0.3}) is in fact a finite
sum and is indeed
 an element of $((\Delta_{3})_{\circ}({\cal A}))(U)$.
It is easy to verify directly that
(\ref{4.0.3}) indeed gives a morphism of ${\cal D}_{X\times
X\times X}$-modules.
Thus we obtain the morphism
we want.

The morphism given by (\ref{4.0.3}) is natural but we need another
expression of the right-hand side of (\ref{4.0.3}) in Section 5.
Note that $\frac{\p}{\p z_{1}}+\frac{\p}{\p z_{2}}+\frac{\p}{\p z_{3}}
\in ((\Delta_{3})_{*}({\cal D}_{X}))(U)$ and its preimage in
${\cal D}_{X}(U_{1}\cap U_{2}\cap U_{3})$ is $\frac{\p}{\p z}$. So the
right-hand side of (\ref{4.0.3}) is equal to
\begin{eqnarray}\label{4.0.4}
\lefteqn{\sum_{i\in {\Bbb N}}
\sum_{k, l\in {\Bbb N}}\sum_{j=0}^{l}
{m_{3}\choose i}(-1)^{j}{l\choose j}\frac{\p^{k}}{\p z^{k}_{1}}
\frac{\p^{j}}{\p z^{j}_{3}}\left(\frac{\p}{\p z_{1}}+
\frac{\p}{\p z_{2}}+\frac{\p}{\p z_{3}}\right)^{l-j}}\nno\\
&&\quad\quad \otimes_{((\Delta_{3})_{*}
({\cal D}_{X}))(U)}(\Delta_{3})_{*}(C^{m_{1}+i,
m_{2}+m_{3}-i}_{kl})\nno\\
&&=\sum_{i\in {\Bbb N}}
\sum_{k, l\in {\Bbb N}}
{m_{3}\choose i}\frac{\p^{k}}{\p z^{k}_{1}}
\left(\frac{\p}{\p z_{1}}+
\frac{\p}{\p z_{2}}\right)^{l}
\otimes_{((\Delta_{3})_{*}
({\cal D}_{X}))(U)}(\Delta_{3})_{*}(C^{m_{1}+i,
m_{2}+m_{3}-i}_{kl})\nno\\
&&=\sum_{i\in {\Bbb N}}
\sum_{k, l\in {\Bbb N}}\sum_{j=0}^{l}
{m_{3}\choose i}{l\choose j}\frac{\p^{k+j}}{\p z^{k+j}_{1}}
\frac{\p^{l-j}}{\p z^{l-j}_{2}}\nno\\
&&\quad\quad \otimes_{((\Delta_{3})_{*}
({\cal D}_{X}))(U)}(\Delta_{3})_{*}(C^{m_{1}+i,
m_{2}+m_{3}-i}_{kl}).
\end{eqnarray}

{}From
\begin{equation}\label{expl-exp1}
\mu((z_{1}-z_{2})^{n+1}(A_{1}\otimes A_{2}))=(z_{1}-z_{2})
\mu((z_{1}-z_{2})^{n}(A_{1}\otimes A_{2}))
\end{equation}
and (\ref{expl-exp}), we obtain
\begin{eqnarray}\label{expl-exp2}
\lefteqn{\sum_{k=0}^{p_{n_{1}+1}}\frac{\p^{k}}{\p z^{k}_{1}}
\otimes_{(\Delta_{*}({\cal
D}_{X}))(U_{1}\times U_{2})}
\Delta_{*}(B^{n_{1}+1}_{k})}\nno\\
&&=(z_{1}-z_{2})\sum_{k=0}^{p_{n_{1}}}\frac{\p^{k}}{\p z^{k}_{1}}
\otimes_{(\Delta_{*}({\cal
D}_{X}))(U_{1}\times U_{2})}
\Delta_{*}(B^{n_{1}}_{k}).
\end{eqnarray}
Since for any $k\in {\Bbb N}$,
$$\frac{\p^{k}}{\p z^{k}_{1}}(z_{1}-z_{2})
\otimes_{(\Delta_{*}({\cal
D}_{X}))(U_{1}\times U_{2})}
\Delta_{*}(B^{n_{1}}_{m})=0,$$
the right-hand side of (\ref{expl-exp2}) is equal to
\begin{equation}\label{expl-exp3}
-\sum_{k=0}^{p_{n_{1}}}k\frac{\p^{k-1}}{\p z^{k-1}_{1}}
\otimes_{(\Delta_{*}({\cal
D}_{X}))(U_{1}\times U_{2})}
\Delta_{*}(B^{n_{1}}_{k}).
\end{equation}
Note that for any $A_{i}\in {\cal A}(U_{1}\cap U_{2})$,
$m_{i}\in {\Bbb N}$, $i=1, \dots, n$,
$$\frac{\p^{m_{i}}}{\p z_{1}^{m_{i}}}\otimes_{(\Delta_{*}({\cal
D}_{X}))(U_{1}\times U_{2})}
\Delta_{*}(A_{i}), \;i=1, \dots, n$$
are linearly independent when
$A_{i}\ne 0$ for $i=1, \dots, n$ and
$m_{i}\ne m_{j}$ for $i\ne j$,
$i, j=1, \dots, n$, and they are equal to $0$ if and only if
$A_{i}=0$, $i=1, \dots, n$.
Thus {}from (\ref{expl-exp2}) and (\ref{expl-exp3}),
we obtain
\begin{equation}\label{expl-exp4}
B^{n_{1}+1}_{k}=-(k+1)B^{n_{1}}_{k+1}
\end{equation}
for any $n_{1}\in \Bbb{Z}$, $k=0, \dots, p_{n_{1}}-1$. Similarly
we have
\begin{eqnarray}\label{expl-exp5}
C^{n_{1}+1, n_{2}+1}_{kl}&=&-(k+1)C^{n_{1}, n_{2}+1}_{k+1,l}\nno\\
&=&-(l+1)C^{n_{1}+1, n_{2}}_{k,l+1}.
\end{eqnarray}
Using (\ref{expl-exp5}) repeatedly, we can write the right-hand side of
(\ref{4.0.4}) as
\begin{eqnarray}\label{4.0.5}
\lefteqn{\sum_{i\in {\Bbb N}}
\sum_{k, l\in {\Bbb N}}\sum_{j=0}^{l}
{m_{3}\choose i}{k+j\choose j}\frac{\p^{k+j}}{\p z^{k+j}_{1}}
\frac{\p^{l-j}}{\p z^{l-j}_{2}}}\nno\\
&&\quad\quad \otimes_{((\Delta_{3})_{*}
({\cal D}_{X}))(U)}(\Delta_{3})_{*}(C^{m_{1}+i-j,
m_{2}+m_{3}-i+j}_{k+j, l-j})\nno\\
&&=\sum_{i'\in {\Bbb N}}
\sum_{k', l'\in {\Bbb N}}
\sum_{\tiny \begin{array}{c}i+k=i'\\ 0\le i, k\le i'\end{array}}
{m_{3}\choose i}{k'\choose k}\frac{\p^{k'}}{\p z^{k'}_{1}}
\frac{\p^{l'}}{\p z^{l'}_{2}}\nno\\
&&\quad\quad \otimes_{((\Delta_{3})_{*}
({\cal D}_{X}))(U)}(\Delta_{3})_{*}(C^{m_{1}+i'-k',
m_{2}+m_{3}-i'+k'}_{k', l'})\nno\\
&&=\sum_{i'\in {\Bbb N}}
\sum_{k', l'\in {\Bbb N}}
{m_{3}+k'\choose i'}\frac{\p^{k'}}{\p z^{k'}_{1}}
\frac{\p^{l'}}{\p z^{l'}_{2}}\nno\\
&&\quad\quad \otimes_{((\Delta_{3})_{*}
({\cal D}_{X}))(U)}(\Delta_{3})_{*}(C^{m_{1}+i'-k',
m_{2}+m_{3}-i'+k'}_{k', l'}),
\end{eqnarray}
where in the last step, we have used the identity
$$\sum_{\tiny \begin{array}{c}i+k=i'\\ 0\le i, k\le i'\end{array}}
{m_{3}\choose i}{k'\choose k}
={m_{3}+k'\choose i'}.$$
By (\ref{4.0.3}), (\ref{4.0.4}) and (\ref{4.0.5}),
we obtain
\begin{eqnarray}\label{operad1}
 \lefteqn{(\mu(\mu(\cdot, \cdot), \cdot))
((z_{1}-z_{2})^{m_{1}}(z_{2}-z_{3})^{m_{2}}(z_{1}-z_{3})^{m_{3}}
(A_{1}\otimes A_{2}\otimes A_{3}))}\nno\\
&&{\displaystyle =\sum_{i\in {\Bbb N}}
\sum_{k, l\in {\Bbb N}}{m_{3}+k\choose i}\frac{\p^{k}}{\p z^{k}_{1}}
\frac{\p^{l}}{\p z^{l}_{2}}\otimes_{((\Delta_{3})_{*}
({\cal D}_{X}))(U)}(\Delta_{3})_{*}(C^{m_{1}+i-k,
m_{2}+m_{3}+k-i}_{kl}).}\nno\\
&&
\end{eqnarray}

We also need another morphism (the other ``composition'')
$$\mu(\cdot, \mu(\cdot, \cdot)):
(j_{3})_{*}j_{3}^{*}({\cal A}\boxtimes {\cal A}\boxtimes {\cal A})
\to (\Delta_{3})_{\circ}{\cal A}$$
of ${\cal D}_{X\times
X\times X}$-modules
defined naturally as follows: Let $U$ be an open subset of $X\times
X\times X$ of the form $U_{1}\times U_{2}\times U_{3}$
where $U_{1}, U_{2}$ and $U_{3}$ are open subsets of $X$.
For any $n_{2}\in {\Bbb Z}$ there exist $s_{n_{2}}\in {\Bbb N}$ and
$$D^{n_{2}}_{l}\in {\cal A}(\Delta^{-1}(U_{2}\times U_{3})) ={\cal
A}(U_{2}\cap U_{3}), \;\;l=0, \dots, s_{n_{2}},$$
such that the element
$\mu((z_{2}-z_{3})^{n_{2}}(A_{1}\otimes A_{2}))
\in (\Delta_{*}({\cal A}))(U_{1}\times U_{2})$
 is equal to
$$\sum_{l=0}^{s_{n_{2}}}\frac{\p^{l}}{\p z^{l}_{2}}\otimes_{(\Delta_{*}({\cal
D}_{X}))(U_{2}\times U_{3})} \Delta_{*}(D^{n_{2}}_{l}).$$
For any $n_{1}\in {\Bbb Z}$, we have
$$\mu((z_{1}-z_{3})^{n_{1}}(A_{1}\otimes D^{n_{2}}_{l}))
=\sum_{k=0}^{r_{n_{1}, n_{2}}}
\frac{\p^{k}}{\p z^{k}_{1}}\otimes_{(\Delta_{*}({\cal
D}_{X}))(U_{1}\times U_{3})} \Delta_{*}(E^{n_{1}n_{2}}_{kl})$$
where $r_{n_{1}, n_{2}}\in {\Bbb N}$ and
$E^{n_{1}n_{2}}_{kl}\in {\cal A}(\Delta^{-1}((U_{1}\cap U_{2})
\times U_{3})
={\cal A}(U_{1}\cap U_{2}\cap U_{3})$,
$k=0, \dots, r_{n_{1}, n_{2}}$, $l=0, \dots, s_{n_{2}}$.
Let $E^{n_{1}n_{2}}_{kl}=0$ for
$k>r_{n_{1}, n_{2}}$ or $l>s_{n_{2}}$.
We define
\begin{eqnarray}\label{operad2}
\lefteqn{(\mu(\cdot, \mu(\cdot, \cdot)))
((z_{1}-z_{2})^{m_{1}}(z_{2}-z_{3})^{m_{2}}(z_{1}-z_{3})^{m_{3}}
(A_{1}\otimes A_{2}\otimes A_{3}))}\nno\\
&&{\displaystyle =\sum_{i\in {\Bbb N}}(-1)^{i}{m_{1}\choose i}
\sum_{k, l\in {\Bbb N}}\frac{\p^{k}}{\p z^{k}_{1}}
\frac{\p^{l}}{\p z^{l}_{2}}\otimes_{((\Delta_{3})_{*}
({\cal D}_{X}))(U)}(\Delta_{3})_{*}(E^{m_{1}+m_{3}-i,
m_{2}+i}_{kl}),}\nno\\
&&
\end{eqnarray}
using the expansion of
$$(z_{1}-z_{2})^{m_{1}}=((z_{1}-z_{3})-(z_{2}-z_{3}))^{m_{1}}$$
in nonnegative powers of $z_{2}-z_{3}$.

We also have an action of the symmetric group $S_{2}$ on
$j_{*}j^{*}({\cal A}\boxtimes {\cal A})$ and an action of $S_{3}$ on
$(j_{3})_{*}j_{3}^{*}({\cal A}\boxtimes {\cal A}\boxtimes {\cal A})$:
In fact, the elements of $S_{2}$ act on $X\times X$, and
we have direct image sheaves of $j_{*}j^{*}({\cal A}\boxtimes {\cal A})$
by these actions. It is easy to see that these direct image sheaves
are the same as $j_{*}j^{*}({\cal A}\boxtimes {\cal A})$, so that
we obtain an action of $S_{2}$ on
$j_{*}j^{*}({\cal A}\boxtimes {\cal A})$. Explicitly, this action can be
described as follows:
For any element
$f(z_{1}, z_{2})(z_{1}-z_{2})^{n}(A_{1}\otimes A_{2})$ of
$(j_{*}j^{*}({\cal A}\boxtimes {\cal A}))(U_{1}\times U_{2})$, we define
$\sigma(f(z_{1}, z_{2})
(z_{1}-z_{2})^{n}(A_{1}\otimes A_{2}))$ to be the element
$f(z_{1}, z_{2})(z_{1}-z_{2})^{n}(A_{2}\otimes A_{1})$
of $(j_{*}j^{*}({\cal A}\boxtimes {\cal A}))(U_{2}\times U_{1})$.
The action of $S_{3}$ on
$(j_{3})_{*}j_{3}^{*}({\cal A}\boxtimes {\cal A}\boxtimes {\cal A})$ is
defined similarly.

\begin{defi}[Beilinson-Drinfeld {\rm \cite{BD}}]\label{chiral}
{\rm Let $X$ be a nonempty open subset of ${\Bbb  C}$.
A {\it chiral algebra over $X$} is
a ${\cal D}_{X}$-module ${\cal A}$ equipped with a morphism
$\mu: j_{*}j^{*}({\cal A}\boxtimes {\cal A})\to \Delta_{\circ}{\cal A}$
of ${\cal D}_{X\times X}$-modules
satisfying the following axioms:

\begin{enumerate}

\item (Skew-symmetry) Let $\sigma_{12}$ be the nontrivial
element of the symmetric group $S_{2}$.
Then
$$
\mu\circ \sigma_{12}=-\mu.
$$

\item (Jacobi identity) Let $\sigma_{12}$ be the element of $S_{3}$ permuting
the first two letters. Then
$$
\mu(\mu(\cdot, \cdot), \cdot)=\mu(\cdot, \mu(\cdot, \cdot))-
\mu(\cdot, \mu(\cdot, \cdot))\circ \sigma_{12}.
$$

\end{enumerate}}
\end{defi}

\renewcommand{\theequation}{\thesection.\arabic{equation}}
\renewcommand{\therema}{\thesection.\arabic{rema}}
\setcounter{equation}{0}
\setcounter{rema}{0}

\section{Vertex algebras without vacuum
and  chiral algebras in the sense of Beilinson-Drin\-feld}

In this section,
 $X$ is a nonempty open subset of ${\Bbb C}$.  We need the following notion:

\begin{defi}\label{va-no-va-x}
{\rm A {\it vertex algebra without vacuum over $X$} is
 a vertex algebra without vacuum $(V, Y, D)$ equipped with an
${\cal O}(X)$-module structure on $V$ such that for any
$f, g\in {\cal O}(X)$ and $u, v\in V$,
$Y(fu, x)gv=fgY(u, x)v$ and $D(fu)=(\frac{\p}{\p z}f)u+fDu$.}
\end{defi}

Any vertex algebra without vacuum tensored with the commutative
associative algebra ${\cal O}(X)$, viewed
as a vertex algebra (see \cite{B}),
is a vertex algebra without vacuum over $X$ (see \cite{B}, \cite{FHL}).

Let $(V, Y, D)$ be a vertex algebra without vacuum over $X$.
For any open (necessarily affine)
subset
$U\subset X$, let
${\cal A}(U)={\cal O}(U)\otimes_{{\cal O}(X)}
 V$. We
define $A(u, U)=1\otimes_{{\cal O}(X)} u\in {\cal A}(U)$.  Note that
${\cal A}(X)={\cal O}(X)\otimes_{{\cal O}(X)}V=V$ and $A(u, X)=u$.
Since ${\cal O}_{X}$ is quasi-coherent,
we see that the ${\cal O}_{X}$-module defined by
$U\mapsto {\cal A}(U)$ for all open subsets $U\subset X$
 is also quasi-coherent.
We define $\frac{\p}{\p z}A(u, U)=A(Du, U)$, so that $\frac{\p}{\p z}$
acts on ${\cal A}(U)$.
 Thus the ${\cal O}_{X}$-module defined by
$U\mapsto {\cal A}(U)$ for all open subsets $U\subset X$ is a
${\cal D}_{X}$-module. We denote this ${\cal
D}_{X}$-module by ${\cal A}$.

For the vertex algebra $V$ without vacuum, we write the
Jacobi identity in the following  form:
\begin{eqnarray*}
&{\displaystyle (x_{1}-x_{2})^{n}Y(u, x_{1}) Y(v, x_{2})-(-x_{2}+x_{1})^{n}
Y(u, x_{2}) Y(v, x_{1})}&\nno\\
&={\displaystyle \res_{x_{0}}x_{0}^{n}x^{-1}_{2}
\delta\left(\frac{x_{1}-x_{0}}{x_{2}}\right)Y(Y(u, x_{0})v, x_{2})}&
\end{eqnarray*}
for $u, v\in V$,  $n\in {\Bbb Z}$.
We rewrite the right-hand side as
\begin{eqnarray*}
\lefteqn{\res_{x_{0}}x_{0}^{n}e^{-x_{0}\frac{\p}{\p x_{1}}}x_{2}^{-1}
\delta\left(\frac{x_{1}}{x_{2}}
\right)Y(Y(u, x_{0})v, x_{2})}\nno\\
&&=\sum_{m\in {\Bbb N}}
\frac{(-1)^{m}}{m!}\frac{\p^{m}}{\p x_{1}^{m}}x^{-1}_{2}
\delta\left(\frac{x_{1}}{x_{2}}\right)Y(u_{m+n}v, x_{2}).
\end{eqnarray*}
We define $\mu: j_{*}j^{*}({\cal A}\boxtimes {\cal A})\to
\Delta_{\circ}({\cal A})$ as follows: For any open subset $U$ of
the form $U_{1}\times U_{2}
\subset X\times X$,
$j_{*}j^{*}({\cal A}\boxtimes {\cal A})(U)$ is spanned by
elements of the form
$$(z_{1}-z_{2})^{n}(A(u, U_{1})\otimes A(v, U_{2})),$$
$u, v\in V$, $n\in {\Bbb Z}$ (or $n<0$).
We define
\begin{eqnarray}\label{mu}
\lefteqn{\mu((z_{1}-z_{2})^{n}(A(u,
U_{1})\otimes A(v, U_{2})))}\nno\\
&&=\sum_{m\in {\Bbb N}}\frac{(-1)^{m}}{m!}
\frac{\p^{m}}{\p z_{1}^{m}}\otimes_{(\Delta_{*}({\cal D}_{X}))(U)}
\Delta_{*}(A(u_{m+n}v, \Delta^{-1}(U))).
\end{eqnarray}
It is easy to verify that $\mu$ is well-defined and is
indeed a morphism of
${\cal D}_{X\times X}$-modules.

We have:

\begin{propo}
The pair $({\cal A}, \mu)$ is a chiral algebra over
$X$ in the sense of Definition \ref{chiral}.
\end{propo}
\pf
We first verify the skew-symmetry using the skew-symmetry for
the vertex algebra without vacuum $V$. By definition, for
$U=U_{1}\times U_{2}$ and $u, v\in V$,
\begin{eqnarray}\label{5.1}
\lefteqn{(\mu\circ \sigma_{12})((z_{1}-z_{2})^{n}(A(u,
U_{1})\otimes A(v, U_{2})))}\nno\\
&&=\mu((z_{1}-z_{2})^{n}(A(v,
U_{2})\otimes A(u, U_{1})))\nno\\
&&=(-1)^{n}\mu((z_{2}-z_{1})^{n}(A(v,
U_{2})\otimes A(u, U_{1})))\nno\\
&&=(-1)^{n}\sum_{m\in {\Bbb N}}\frac{(-1)^{m}}{m!}
\frac{\p^{m}}{\p z_{2}^{m}}\otimes_{(\Delta_{*}({\cal D}_{X}))(U)}
\Delta_{*}(A(v_{m+n}u, \Delta^{-1}(\sigma_{12}(U)))
\nno\\
&&=\sum_{m\in {\Bbb N}}\frac{(-1)^{m+n}}{m!}
\frac{\p^{m}}{\p z_{2}^{m}}\otimes_{(\Delta_{*}({\cal D}_{X}))(U)}
\Delta_{*}(A(v_{m+n}u,  \Delta^{-1}(U)))\nno\\
&&=\sum_{m\in {\Bbb N}}\frac{(-1)^{m+n}}{m!}
\left(-\frac{\p}{\p z_{1}}+\left(\frac{\p}{\p z_{1}}
+\frac{\p}{\p z_{2}}\right)
\right)^{m}\nno\\
&&\hspace{2em}\otimes_{(\Delta_{*}({\cal D}_{X}))(U)}
\Delta_{*}(A(v_{m+n}u,  \Delta^{-1}(U)))\nno\\
&&=\sum_{m\in {\Bbb N}}\sum_{k\in {\Bbb N}}
\frac{(-1)^{m+n}}{m!}{m\choose k}
\left(-\frac{\p}{\p z_{1}}\right)^{m-k}
\left(\frac{\p}{\p z_{1}}+\frac{\p}{\p z_{2}}
\right)^{k}\nno\\
&&\hspace{2em}\otimes_{(\Delta_{*}({\cal D}_{X}))(U)}
\Delta_{*}(A(v_{m+n}u,  \Delta^{-1}(U)))\nno\\
&&=\sum_{l\in {\Bbb N}}\sum_{k\in {\Bbb N}}
\frac{(-1)^{k+n}}{l!k!}
\left(\frac{\p}{\p z_{1}}\right)^{l}
\left(\frac{\p}{\p z_{1}}+\frac{\p}{\p z_{2}}
\right)^{k}\nno\\
&&\hspace{2em}\otimes_{(\Delta_{*}({\cal D}_{X}))(U)}
\Delta_{*}(A(v_{k+l+n}u,  \Delta^{-1}(U))).
\end{eqnarray}
Since $\frac{\p}{\p z_{1}}+\frac{\p}{\p z_{2}}\in
(\Delta_{*}({\cal D}_{X}))(U)$ and its preimage in
$${\cal D}_{X}(\Delta^{-1}(U))
={\cal D}_{X}(U_{1}\cap U_{2})$$
is
$\frac{\p}{\p z}$, the right-hand side of (\ref{5.1}) is equal to
\begin{eqnarray}\label{5.2}
\lefteqn{\sum_{l\in {\Bbb N}}\sum_{k\in {\Bbb N}}
\frac{(-1)^{k+n}}{l!k!}
\left(\frac{\p}{\p z_{1}}\right)^{l}
\otimes_{(\Delta_{*}({\cal D}_{X}))(U)}
\Delta_{*}\left(\left(\frac{\p}{\p z}\right)^{k}
A(v_{k+l+n}u,  \Delta^{-1}(U))\right)}\nno\\
&&=\sum_{l\in {\Bbb N}}\sum_{k\in {\Bbb N}}
\frac{(-1)^{k+n}}{l!k!}
\left(\frac{\p}{\p z_{1}}\right)^{l}
\otimes_{(\Delta_{*}({\cal D}_{X}))(U)}
\Delta_{*}\left(
A(D^{k}v_{k+l+n}u,  \Delta^{-1}(U))\right)\nno\\
&&=-\sum_{l\in {\Bbb N}}
\frac{(-1)^{l}}{l!}
\left(\frac{\p}{\p z_{1}}\right)^{l}\nno\\
&&\hspace{2em}\otimes_{(\Delta_{*}({\cal D}_{X}))(U)}
\Delta_{*}\left(
A\left(\sum_{k\in {\Bbb N}}\frac{(-1)^{k+l+n+1}}{k!}D^{k}v_{k+l+n}u,
 \Delta^{-1}(U)\right)\right).\nno\\
&&
\end{eqnarray}
But {}from the skew-symmetry (\ref{skew}),
we obtain its component form
$$\sum_{k\in {\Bbb N}}\frac{(-1)^{k+m+1}}{k!}D^{k}v_{k+m}u
=u_{m}v$$
for all $m\in \Bbb{Z}$. Thus the right-hand side of (\ref{5.2}) is
equal to
\begin{equation}\label{5.3}
-\sum_{l\in {\Bbb N}}
\frac{(-1)^{l}}{l!}
\left(\frac{\p}{\p z_{1}}\right)^{l}
\otimes_{(\Delta_{*}({\cal D}_{X}))(U)}
\Delta_{*}\left(
A\left(u_{l+n}v,
 \Delta^{-1}(U)\right)\right).
\end{equation}

On the other hand,
\begin{eqnarray}\label{5.4}
\lefteqn{(-\mu)((z_{1}-z_{2})^{n}A(u,
U_{1})\otimes A(v, U_{2}))}\nno\\
&&=-\mu((z_{1}-z_{2})^{n}A(u,
U_{1})\otimes A(v, U_{2}))\nno\\
&&=-\sum_{m\in {\Bbb N}}\frac{(-1)^{m}}{m!}
\frac{\p^{m}}{\p z_{1}^{m}}\otimes_{(\Delta_{*}({\cal D}_{X}))(U)}
\Delta_{*}(A(u_{m+n}v,  \Delta^{-1}(U))).
\end{eqnarray}
We see that the right-hand side of (\ref{5.4}) is equal to
(\ref{5.3}). By the calculation {}from (\ref{5.1}) to (\ref{5.4}), we see
that the left-hand side of (5.1) and the left-hand side of (5.4) are equal,
proving the skew-symmetry.

Next we prove the Jacobi identity using the Jacobi identity for
vertex algebra without vacuum. By (\ref{mu}) and
(\ref{operad2}), for any open subset
$U$ of the form $U_{1}\times U_{2}\times U_{3}
\subset X\times X\times X$,
\begin{eqnarray}\label{j1}
\lefteqn{(\mu(\cdot, \mu(\cdot, \cdot)))
((z_{1}-z_{2})^{m_{1}}(z_{2}-z_{3})^{m_{2}}
(z_{1}-z_{3})^{m_{3}}\cdot}\nno\\
&&\hspace{8em}\cdot (A(u, U_{1})\otimes A(v, U_{2})
\otimes A(w, U_{3})))\nno\\
&&=\sum_{i\in {\Bbb N}}(-1)^{i}{m_{1}\choose i}
\sum_{k, l\in {\Bbb N}}\frac{(-1)^{k+l}}{k!l!}\frac{\p^{k}}{\p z^{k}_{1}}
\frac{\p^{l}}{\p z^{l}_{2}}\nno\\
&&\hspace{3em}\otimes_{((\Delta_{3})_{*}
({\cal D}_{X}))(U)}(\Delta_{3})_{*}(A(u_{m_{1}+m_{3}-i+k}v_{
m_{2}+i+l}w, \Delta_{3}^{-1}(U)))\nno\\
&&=
\sum_{k, l\in {\Bbb N}}\frac{(-1)^{k+l}}{k!l!}\frac{\p^{k}}{\p z^{k}_{1}}
\frac{\p^{l}}{\p z^{l}_{2}}\nno\\
&&\hspace{3em}\otimes_{((\Delta_{3})_{*}
({\cal D}_{X}))(U)}(\Delta_{3})_{*}\biggl(A\biggl(
\sum_{i\in {\Bbb N}}(-1)^{i}{m_{1}\choose i}\cdot\nno\\
&&\hspace{13em}\cdot u_{m_{1}+m_{3}-i+k}v_{
m_{2}+i+l}w, \Delta_{3}^{-1}(U)\biggr)\biggr).\nno\\
&&
\end{eqnarray}
Similarly, we have
\begin{eqnarray}\label{j2}
\lefteqn{(\mu(\cdot, \mu(\cdot, \cdot))\circ\sigma_{12})
((z_{1}-z_{2})^{m_{1}}(z_{2}-z_{3})^{m_{2}}
(z_{1}-z_{3})^{m_{3}}\cdot}\nno\\
&&\hspace{8em}\cdot(A(u, U_{1})\otimes A(v, U_{2})
\otimes A(w, U_{3})))\nno\\
&&=
\sum_{k, l\in {\Bbb N}}\frac{(-1)^{k+l}}{k!l!}\frac{\p^{k}}{\p z^{k}_{1}}
\frac{\p^{l}}{\p z^{l}_{2}}\nno\\
&&\hspace{3em}\otimes_{((\Delta_{3})_{*}
({\cal D}_{X}))(U)}(\Delta_{3})_{*}\biggl(A\biggl((-1)^{m_{1}}
\sum_{i\in {\Bbb N}}(-1)^{i}{m_{1}\choose i}\cdot\nno\\
&&\hspace{13em}\cdot v_{m_{1}+m_{2}-i+l}u_{
m_{3}+i+k}w, \Delta_{3}^{-1}(U)\biggr)\biggr).\nno\\
&&
\end{eqnarray}
On the other hand, by (\ref{mu}) and (\ref{operad1}), we have
\begin{eqnarray}\label{j3}
\lefteqn{(\mu(\mu(\cdot, \cdot), \cdot)))
((z_{1}-z_{2})^{m_{1}}(z_{2}-z_{3})^{m_{2}}
(z_{1}-z_{3})^{m_{3}}\cdot}\nno\\
&&\hspace{4em}\cdot(A(u, U_{1})\otimes A(v, U_{2})
\otimes A(w, U_{3})))\nno\\
&&=\sum_{k, l\in {\Bbb N}}\frac{(-1)^{k+l}}{k!l!}\frac{\p^{k}}{\p z^{k}_{1}}
\frac{\p^{l}}{\p z^{l}_{2}}\nno\\
&&\hspace{3em}\otimes_{((\Delta_{3})_{*}
({\cal D}_{X}))(U)}(\Delta_{3})_{*}\biggl(A\biggl(
\sum_{i\in {\Bbb N}}{m_{3}+k\choose i}\cdot\nno\\
&&\hspace{8em} \cdot(u_{m_{1}+i}v)_{
m_{3}+m_{2}-i+k+l}w, \Delta_{3}^{-1}(U)\biggr)\biggr).
\end{eqnarray}
{}From the Jacobi identity for vertex algebras without vacuum, we can
obtain its component form
\begin{eqnarray}\label{jac-comp}
\lefteqn{\sum_{i\in {\Bbb N}}{m\choose i}(u_{l+i}v)_{
m+n-i}w=}\nno\\
&&=\sum_{i\in {\Bbb N}}(-1)^{i}{l\choose i}u_{m+l-i}v_{
n+i}w\nno\\
&&\quad -(-1)^{l}
\sum_{i\in {\Bbb N}}(-1)^{i}{l\choose i}v_{n+l-i}u_{
m+i}w
\end{eqnarray}
for all $l, m, n\in {\Bbb Z}$ and $u, v, w\in V$
(see \cite{FLM}). {}From (\ref{jac-comp}),
we obtain
\begin{eqnarray}\label{component}
\lefteqn{\sum_{i\in {\Bbb N}}{m_{3}+k\choose i}(u_{m_{1}+i}v)_{
m_{3}+m_{2}-i+k+l}w=}\nno\\
&&=\sum_{i\in {\Bbb N}}(-1)^{i}{m_{1}\choose i}u_{m_{1}+m_{3}-i+k}v_{
m_{2}+i+l}w\nno\\
&&\quad -(-1)^{l}
\sum_{i\in {\Bbb N}}(-1)^{i}{m_{1}\choose i}v_{m_{1}+m_{2}-i+l}u_{
m_{3}+i+k}w
\end{eqnarray}
for any $m_{1}, m_{2}, m_{3}\in {\Bbb Z}$, $k, l\in {\Bbb N}$ and
$u, v, w\in V$. By (\ref{j1})--(\ref{component}), we obtain
\begin{eqnarray*}
\lefteqn{(\mu(\mu(\cdot, \cdot), \cdot)))
((z_{1}-z_{2})^{m_{1}}(z_{2}-z_{3})^{m_{2}}
(z_{1}-z_{3})^{m_{3}}\cdot}\nno\\
&&\hspace{4em}\cdot(A(u, U_{1})\otimes A(v, U_{2})
\otimes A(w, U_{3})))\nno\\
&&=(\mu(\cdot, \mu(\cdot, \cdot)))
((z_{1}-z_{2})^{m_{1}}(z_{2}-z_{3})^{m_{2}}
(z_{1}-z_{3})^{m_{3}}\cdot\nno\\
&&\hspace{3em}\cdot(A(u, U_{1})\otimes A(v, U_{2})
\otimes A(w, U_{3})))\nno\\
&&\quad -(\mu(\cdot, \mu(\cdot, \cdot))\circ\sigma_{12})
((z_{1}-z_{2})^{m_{1}}(z_{2}-z_{3})^{m_{2}}
(z_{1}-z_{3})^{m_{3}}\cdot\nno\\
&&\hspace{3em}\cdot(A(u, U_{1})\otimes A(v, U_{2})
\otimes A(w, U_{3}))),
\end{eqnarray*}
proving the Jacobi identity for ${\cal A}$.
\epfv

Conversely, suppose that we have a chiral algebra $({\cal A}, \mu)$ over
a nonempty open subset $X$ of ${\Bbb C}$
in the sense of Definition \ref{chiral}, and let $V={\cal A}(X)$
be the space of all global sections over $X$.
Then $V$ is an ${\cal O}(X)$-module and $\frac{\p}{\p z}$ acts on $V$.
For any element $v\in V$, we define
$Dv=\frac{\p}{\p z}v$. For any $f\in {\cal O}(X)$ and $u\in V$, we have
$D(fu)=(\frac{\p}{\p z}f)u+fDu$. For $u, v\in V$,
there exist $B^{n}_{m}(u, v)\in
V$, $n\in {\Bbb Z}$, $m=1, \dots, p_{n}$, such
that
\begin{equation}\label{expl-exp0}
\mu((z_{1}-z_{2})^{n}(u\otimes v))
=\sum_{m=0}^{p_{n}}\frac{\p^{m}}{\p z^{m}_{1}}\otimes_{(\Delta_{*}({\cal
D}_{X}))(X\times X)}
\Delta_{*}(B^{n}_{m}(u, v)),
\end{equation}
and the nonzero $B^{n}_{m}(u, v)$ are uniquely determined.
The equation (\ref{expl-exp4}) in this case becomes
\begin{equation}\label{expl-exp6}
B^{n+1}_{m}(u, v)=-(m+1)B^{n}_{m+1}(u, v)
\end{equation}
for any $u, v\in V$, $n\in {\Bbb Z}$ and $m=0, \dots, p_{n}-1$.

We define
$u_{n}v=B^{n}_{0}(u, v)$ and
$$Y(u, x)v=\sum_{n\in {\Bbb Z}}B^{n}_{0}(u, v)x^{-n-1}.$$
{}From the definition, we have
$Y(fu, x)gv=fgY(u, x)v$ for any $f, g\in {\cal O}(X)$ and $u, v\in V$.

\begin{propo}
The triple $(V, Y, D)$
is a vertex algebra without vacuum over $X$.
\end{propo}
\pf
When $n>p_{0}$, we have
\begin{eqnarray*}
\lefteqn{\mu((z_{1}-z_{2})^{n}(u\otimes v))=}\nno\\
&&=(z_{1}-z_{2})^{n}\mu(u\otimes v)\nno\\
&&=(z_{1}-z_{2})^{n}\sum_{m=0}^{p_{0}}
\frac{\p^{m}}{\p z^{m}_{1}}\otimes_{(\Delta_{*}({\cal
D}_{X}))(X\times X)}
\Delta_{*}(B^{0}_{m}(u, v))\nno\\
&&=0.
\end{eqnarray*}
Thus for $n>p_{0}$, $u_{n}v=B^{n}_{0}(u, v)=0.$

To show the $D$-derivative property, we need only show its
component form:
\begin{equation}\label{l-1-0}
-(n+1)B^{n}_{0}(u, v)=B^{n+1}_{0}(Du, v)
\end{equation}
for all $n\in {\Bbb Z}$. In fact,
for any $n\in {\Bbb Z}$,
\begin{eqnarray}\label{l-1-1}
\lefteqn{\frac{\p}{\p z_{1}}\mu((z_{1}-z_{2})^{n+1}(u\otimes v))=}\nno\\
&&=(n+1)\mu((z_{1}-z_{2})^{n}(u\otimes v))
+\mu((z_{1}-z_{2})^{n+1}\frac{\p}{\p z_{1}}(u\otimes v))\nno\\
&&=(n+1)\mu((z_{1}-z_{2})^{n}(u\otimes v))
+\mu\left((z_{1}-z_{2})^{n+1}\left(\left(\frac{\p}{\p z}u\right)\otimes v
\right)\right)\nno\\
&&=(n+1)\mu((z_{1}-z_{2})^{n}(u\otimes v))
+\mu((z_{1}-z_{2})^{n+1}((Du)\otimes v)).
\end{eqnarray}
On the other hand,
\begin{eqnarray}\label{l-1-2}
\lefteqn{\frac{\p}{\p z_{1}}\mu((z_{1}-z_{2})^{n+1}(u\otimes v))=}\nno\\
&&= \frac{\p}{\p z_{1}}\left(
\sum_{m=0}^{p_{n}}\frac{\p^{m}}{\p z^{m}_{1}}\otimes_{(\Delta_{*}({\cal
D}_{X}))(X\times X)}
\Delta_{*}(B^{n}_{m}(u, v))\right)
\nno\\
&&=\sum_{m=0}^{p_{n}}\frac{\p^{m+1}}{\p z^{m+1}_{1}}
\otimes_{(\Delta_{*}({\cal
D}_{X}))(X\times X)}
\Delta_{*}(B^{n}_{m}(u, v)).
\end{eqnarray}
{}From (\ref{l-1-1}), (\ref{l-1-2}) and the expansions of
$\mu((z_{1}-z_{2})^{n}(u\otimes v))$ and
$\mu((z_{1}-z_{2})^{n+1}((Du)\otimes v))$, we obtain
\begin{eqnarray}\label{l-1-3}
\lefteqn{(n+1)\sum_{m=0}^{p_{n}}\frac{\p^{m}}{\p z^{m}_{1}}
\otimes_{(\Delta_{*}({\cal
D}_{X}))(X\times X)}
\Delta_{*}(B^{n}_{m}(u, v))}\nno\\
&&\quad +\sum_{m=0}^{p_{n+1}}\frac{\p^{m}}{\p z^{m}_{1}}
\otimes_{(\Delta_{*}({\cal
D}_{X}))(X\times X)}
\Delta_{*}(B^{n+1}_{m}(u, v))\nno\\
&&=(n+1)\mu((z_{1}-z_{2})^{n}(u\otimes v))
+\mu((z_{1}-z_{2})^{n+1}((Du)\otimes v))\nno\\
&&=\sum_{m=0}^{p_{n}}\frac{\p^{m+1}}{\p z^{m+1}_{1}}
\otimes_{(\Delta_{*}({\cal
D}_{X}))(X\times X)}
\Delta_{*}(B^{n}_{m}(u, v)).
\end{eqnarray}
The equality
(\ref{l-1-3}) implies
$$(n+1)B^{n}_{0}(u, v)+B^{n+1}_{0}(u, v)=0$$
which is equivalent to (\ref{l-1-0}).

We now prove that
the skew-symmetry for vertex algebras without vacuum
is satisfied by $V$,
using the skew-symmetry for $({\cal A}, \mu)$.
By (\ref{expl-exp0}) and the skew-symmetry for $({\cal A}, \mu)$, we obtain
\begin{eqnarray}\label{skew1}
\lefteqn{\sum_{m=0}^{p_{n}}\frac{\p^{m}}{\p z^{m}_{1}}
\otimes_{(\Delta_{*}({\cal
D}_{X}))(X\times X)}
\Delta_{*}(B^{n}_{m}(u, v))}\nno\\
&&=-(-1)^{n}\sum_{m=0}^{q_{n}}\frac{\p^{m}}{\p z^{m}_{2}}
\otimes_{(\Delta_{*}({\cal
D}_{X}))(X\times X)}
\Delta_{*}(B^{n}_{m}(v, u)).
\end{eqnarray}
The right-hand side of (\ref{skew1}) is equal to
\begin{eqnarray}\label{skew2}
\lefteqn{-(-1)^{n}\sum_{m=0}^{q_{n}}\left(-\frac{\p}{\p z_{1}}+
\left(\frac{\p}{\p z_{1}}+\frac{\p}{\p z_{2}}\right)\right)^{m}
\otimes_{(\Delta_{*}({\cal
D}_{X}))(X\times X)}
\Delta_{*}(B^{n}_{m}(v, u))}\nno\\
&&=(-1)^{n+1}\sum_{m=0}^{q_{n}}\sum_{k=0}^{m}
(-1)^{m-k}{m \choose k}\cdot\nno\\
&&\hspace{3em}\cdot \frac{\p^{m-k}}{\p z^{m-k}_{1}}
\left(\frac{\p}{\p z_{1}}+\frac{\p}{\p z_{2}}\right)^{k}
\otimes_{(\Delta_{*}({\cal
D}_{X}))(X\times X)}
\Delta_{*}(B^{n}_{m}(v, u))\nno\\
&&=(-1)^{n+1}\sum_{m=0}^{q_{n}}\sum_{k=0}^{m}
(-1)^{m-k}{m \choose k}\frac{\p^{m-k}}{\p z^{m-k}_{1}}
\otimes_{(\Delta_{*}({\cal
D}_{X}))(X\times X)}
\Delta_{*}\left(\frac{d^{k}}{d z^{k}}B^{n}_{m}(v, u)\right)\nno\\
&&=(-1)^{n+1}\sum_{m=0}^{q_{n}}\sum_{k=0}^{m}
(-1)^{m-k}{m \choose k}\frac{\p^{m-k}}{\p z^{m-k}_{1}}
\otimes_{(\Delta_{*}({\cal
D}_{X}))(X\times X)}
\Delta_{*}(D^{k}B^{n}_{m}(v, u)).\nno\\
&&
\end{eqnarray}
Comparing the left-hand side of (\ref{skew1}) and the right-hand side of
(\ref{skew2}), the terms involving $\frac{\p^{0}}{\p z_{1}^{0}}$ give
\begin{equation}\label{skew3}
B^{n}_{0}(u, v)=\sum_{k=0}^{q_{n}}(-1)^{n+1}D^{k}B^{n}_{k}(v, u).
\end{equation}
Using (\ref{expl-exp6}) repeatedly, we obtain
\begin{equation}\label{skew4}
B^{n}_{k}(v, u)=\frac{(-1)^{k}}{k!}B_{0}^{k+n}(v, u).
\end{equation}
Combining (\ref{skew3}) and (\ref{skew4}), we obtain
\begin{equation}\label{skew5}
B^{n}_{0}(u, v)=\sum_{k=0}^{q_{n}}\frac{(-1)^{k+n+1}}{k!}
D^{k}B^{k+n}_{0}(v, u).
\end{equation}
Since when $k> q_{n}$, we have
\begin{eqnarray*}
\lefteqn{\mu((z_{1}-z_{2})^{k+n}(v\otimes u))=}\nno\\
&&=(z_{1}-z_{2})^{k}\mu((z_{1}-z_{2})^{n}(v\otimes u))\nno\\
&&=(z_{1}-z_{2})^{k}\sum_{m=0}^{q_{n}}\frac{\p^{m}}{\p z_{2}^{m}}
\otimes_{(\Delta_{*}({\cal
D}_{X}))(X\times X)}
\Delta_{*}(D^{k}B^{n}_{m}(v, u))\nno\\
&&=0.
\end{eqnarray*}
Thus $B^{k+n}_{0}(v, u)=0$ when $k>q_{n}$. So (\ref{skew5}) becomes
$$B^{n}_{0}(u, v)=\sum_{k\in {\Bbb N}}
\frac{(-1)^{k+n+1}}{k!}D^{k}B^{k+n}_{0}(v, u),$$
or equivalently
$$u_{n}v=\sum_{k\in {\Bbb N}}
\frac{(-1)^{k+n+1}}{k!}D^{k}v_{k+n}u$$
which is the component form of the skew-symmetry.

Finally, we prove that the Jacobi identity for vertex algebra
without vacuum is satisfied by $V$, using the
Jacobi identity for $({\cal A}, \mu)$.
By (\ref{operad1}), (\ref{operad2}) and the Jacobi identity for
 $({\cal A}, \mu)$, we obtain the following
identity:
\begin{eqnarray}\label{comp-jac}
\lefteqn{\sum_{i\in {\Bbb N}}
\sum_{k, l\in {\Bbb N}}{m_{3}+k\choose i}\frac{\p^{k}}{\p z^{k}_{1}}
\frac{\p^{l}}{\p z^{l}_{2}}}\nno\\
&&\quad \otimes_{((\Delta_{3})_{*}
({\cal D}_{X}))(U)}(\Delta_{3})_{*}(B^{m_{2}+m_{3}-k-i}_{l}
(B_{k}^{m_{1}+i+k}(u,
v), w))\nno\\
&&=\sum_{i\in {\Bbb N}}(-1)^{i}{m_{1}\choose i}
\sum_{k, l\in {\Bbb N}}\frac{\p^{k}}{\p z^{k}_{1}}
\frac{\p^{l}}{\p z^{l}_{2}}\otimes_{((\Delta_{3})_{*}
({\cal D}_{X}))(U)}\nno\\
&&\hspace{3em}\otimes_{((\Delta_{3})_{*}
({\cal D}_{X}))(U)}(\Delta_{3})_{*}(B^{m_{1}+m_{3}-i}_{k}(u,
B^{m_{2}+i}_{l}(v, w)))\nno\\
&&\quad -(-1)^{m_{1}}\sum_{i\in {\Bbb N}}(-1)^{i}{m_{1}\choose i}
\sum_{k, l\in {\Bbb N}}\frac{\p^{k}}{\p z^{k}_{1}}
\frac{\p^{l}}{\p z^{l}_{2}}\otimes_{((\Delta_{3})_{*}
({\cal D}_{X}))(U)}\nno\\
&&\hspace{3em}\otimes_{((\Delta_{3})_{*}
({\cal D}_{X}))(U)}(\Delta_{3})_{*}(B^{m_{1}+m_{2}-i}_{l}(v,
B^{m_{3}+i}_{k}(u, w))).
\end{eqnarray}
This equality implies
 \begin{eqnarray*}
\lefteqn{\sum_{i\in {\Bbb N}}{m_{3}\choose i}
B^{m_{2}+m_{3}-i}_{0}(B_{0}^{m_{1}+i}(u, v), w)=}\nno\\
&&=\sum_{i\in {\Bbb N}}(-1)^{i}{m_{1}\choose i}B^{m_{1}+m_{3}-i}_{0}(u,
B^{m_{2}+i}_{0}(v,
w))\nno\\
&&\quad \quad -(-1)^{m_{1}}\sum_{i\in {\Bbb N}}(-1)^{i}{m_{1}\choose i}
B^{m_{1}+m_{2}-i}_{0}(v, B^{m_{3}+i}_{0}(u,
w)))
\end{eqnarray*}
or equivalently, the component form of the Jacobi identity
for vertex algebra without vacuum:
\begin{eqnarray*}
\lefteqn{\sum_{i\in {\Bbb N}}{m_{3}\choose i}
(u_{m_{1}+i}v)_{m_{2}+m_{3}-i}w=}\nno\\
&&=\sum_{i\in {\Bbb N}}(-1)^{i}{m_{1}\choose i}u_{m_{1}+m_{3}-i}
v_{m_{2}+i}w\nno\\
&&\quad \quad -(-1)^{m_{1}}\sum_{i\in {\Bbb N}}(-1)^{i}{m_{1}\choose i}
v_{m_{1}+m_{2}-i}u_{m_{3}+i}w.
\end{eqnarray*}

Since
 $(V, Y, D)$ satisfies
the lower-truncation condition for
vertex operators, the $D$-derivative property, the skew-symmetry
and the Jacobi
identity for vertex algebras without vacuum,
it is a vertex algebra without vacuum.
We already know that $V$ is an ${\cal O}(X)$-module
and for any $f, g\in {\cal O}(X)$,
$u, v\in V$, $Y(fu, x)gv=fgY(u, x)v$ and $D(fu)=(\frac{\p}{\p z}f)u+fDu$. So
$(V, Y, D)$ is a vertex algebra without vacuum over $X$.
\epfv

{}From the  two propositions above, we obtain the following
main result of this paper:

\begin{theo}
Let $X$ be a nonempty open subset of ${\Bbb C}$.
The category of vertex algebras without vacuum over $X$ and the category of
chiral algebras over $X$ are equivalent.\epf
\end{theo}

 {\small \sc  Department of Mathematics, Rutgers University,
New Brunswick, NJ 08903}

{\em E-mail address}: yzhuang@math.rutgers.edu

\vskip 1em

{\small \sc Department of Mathematics, Rutgers University,
New Brunswick, NJ 08903}

{\em E-mail address}: lepowsky@math.rutgers.edu

\end{document}